\documentclass[12pt]{article}

\advance\voffset by -0.25cm
\usepackage{graphicx,epsfig,amsmath,amsfonts,amsthm,setspace,url,multirow}

\newtheorem{theorem}{Theorem}[section]

\newtheorem{proposition}[theorem]{Proposition}

\theoremstyle{definition}
\newtheorem{definition}[theorem]{Definition}
\newtheorem{example}[theorem]{Example}

\begin{document}

\newcommand{\real}{{\ensuremath{\mathbb{R}}}}
\newcommand{\qq}{{\ensuremath{\mathbb{Q}}}}
\newcommand{\one}{{\ensuremath{\mathbb{I}}}}

\newcommand{\cA}{{\ensuremath{\cal A}}}
\newcommand{\cB}{{\ensuremath{\cal B}}}
\newcommand{\cC}{{\ensuremath{\cal C}}}
\newcommand{\cD}{{\ensuremath{\cal D}}}
\newcommand{\cF}{{\ensuremath{\cal F}}}
\newcommand{\cG}{{\ensuremath{\cal G}}}
\newcommand{\cL}{{\ensuremath{\cal L}}}
\newcommand{\cN}{{\ensuremath{\cal N}}}
\newcommand{\cP}{{\ensuremath{\cal P}}}

\newcommand{\myE}{{\ensuremath{\mathbb{E}}}}
\newcommand{\myP}{{\ensuremath{\mathbb{P}}}}
\newcommand{\myQ}{{\ensuremath{\mathbb{Q}}}}

\newcommand{\myI}{{\ensuremath{\rm I}}}

\newcommand{\dd}{{\ensuremath{\rm d}}}

\newcommand{\hsp}{\hspace{0.2mm}}
\newcommand{\bysame}{\raisebox{1.5mm}{\underline{\hspace{2.75em}}}}

\newcommand{\var}{\mbox{\rm var}}
\newcommand{\cov}{\mbox{\rm cov}}

\newcommand{\done}{\hfill $\Box$}

\newcommand{\marginal}[1]{\marginpar{\raggedright\scriptsize #1}}

\newcommand{\mypath}{}

\begin{center}
{\LARGE \bf

Combining Predictive Distributions

}

\bigskip
{\bf Tilmann Gneiting, University of Heidelberg} \\
{\bf Heidelberg, Germany} 

\medskip
{\bf Roopesh Ranjan, GE Research, Bangalore, India} 

\bigskip
{\bf \today}
\end{center}

\medskip
\begin{abstract}

Predictive distributions need to be aggregated when probabilistic
forecasts are merged, or when expert opinions expressed in terms of
probability distributions are fused.  We take a prediction space
approach that applies to discrete, mixed discrete-continuous and
continuous predictive distributions alike, and study combination
formulas for cumulative distribution functions from the perspectives
of coherence, probabilistic and conditional calibration, and
dispersion.  Both linear and non-linear aggregation methods are
investigated, including generalized, spread-adjusted and
beta-transformed linear pools.  The effects and techniques are
demonstrated theoretically, in simulation examples, and in case
studies on density forecasts for S\&P 500 returns and daily maximum
temperature at Seattle-Tacoma Airport.

\end{abstract}

\smallskip
\section{Introduction} \label{sec:intro}

Probabilistic forecasts aim to provide calibrated and sharp predictive
distributions for future quantities or events of interest.  As they
admit the assessment of forecast uncertainty and allow for optimal
decision making, probabilistic forecasts continue to gain prominence
in a wealth of applications, ranging from economics and finance to
meteorology and climatology (Gneiting 2008).  The general goal is to
maximize the sharpness of the predictive distributions subject to
calibration (Murphy and Winkler 1987; Gneiting, Balabdaoui and Raftery
2007).  For a real-valued outcome, a probabilistic forecast can be
represented in the form of a predictive cumulative distribution
function, which might be discrete, mixed discrete-continuous or
continuous, with the latter case corresponding to density forecasts.

In many situations, complementary or competing probabilistic forecasts
from dependent or independent information sources are available.  For
example, the individual forecasts might stem from distinct experts,
organizations or statistical models.  The prevalent method for
aggregating the individual predictive distributions into a single
combined forecast is the linear pool (Stone 1961).  While other
methods for combining predictive distributions are available (Genest
and Zidek 1986; Clemen and Winkler 1999; 2007), the linear pool is
typically the method of choice, with the pioneering work of Winkler
(1968) and Zarnowitz (1969), and recent papers by Mitchell and Hall
(2005), Wallis (2005), Hall and Mitchell (2007), Jore, Mitchell and
Vahey (2010), Kascha and Ravazzolo (2010) and Garratt et al.~(2011)
being examples in the case of density forecasts.  Similarly, linear
pools have been applied successfully to combine discrete predictive
distributions; for recent reviews, see Ranjan and Gneiting (2010) and
Clements and Harvey (2011).

Despite the ubiquitous success of the linear pool in a vast number of
applications, fragmented recent work points at potential shortcomings
and limitations.  Hora (2004) proved that any nontrivial convex
combination of two calibrated density forecasts is uncalibrated.  In
the discrete case, Dawid, DeGroot and Mortera (1995) and Ranjan and
Gneiting (2010) showed that linear combination formulas with strictly
positive coefficients fail to be coherent and demonstrated potential
improvement under nonlinear aggregation.

Our goal here is to unify these results and to extend them in various
directions.  Towards this end, we develop novel theoretical approaches
to studying combination formulas and aggregation methods, where we
think of an aggregation method as a class of combination formulas.
For example, the traditional linear pools comprises the linear
combination formulas with nonnegative weights that sum to one.
Technically, we operate in terms of cumulative distribution functions,
which allows us, in contrast to earlier work, to provide a unified
treatment of all real-valued predictands, including the cases of
density forecasts, mixed discrete-continuous predictive distributions,
probability mass functions for count data and probability forecasts of
a dichotomous event, all of which are important in applications.  The
extant literature compares combination formulas by examining whether
or not they possess certain analytic characteristics, such as the
strong setwise function and external Bayes properties (Genest and
Zidek 1986; French and R\'{i}os Insua 2000).  Again in contrast to
much of the earlier work, we assess combination formulas and
aggregation methods from the perspectives of coherence, calibration
and dispersion.

Section \ref{sec:theory} sets the stage by introducing the key tool of
a prediction space, which is a probability space tailored to the study
of forecasts and combination formulas.  We revisit the work of
Gneiting et al.~(2007) and Ranjan and Gneiting (2010) in the
prediction space setting and show, perhaps surprisingly, that if the
outcome is binary, probabilistic calibration and conditional
calibration are equivalent.  Throughout the technical parts of the
paper, information sets are represented by $\sigma$-algebras, and a
predictive distribution is considered ideal if it agrees with the
conditional distribution of the predictand, given the information
basis or $\sigma$-algebra at hand.
 
Section \ref{sec:methods} is devoted to the study of specific, linear
and non-linear combination formulas and aggregation methods.  A key
result is, roughly, that dispersion tends to increase under linear
pooling.  This helps explain the success of linear combination
formulas in aggregating underdispersed component distributions, and
leads us to unify and strengthen the aforementioned results of Dawid
et al.~(1995), Hora (2004) and Ranjan and Gneiting (2010).  In
particular, any linear combination formula with strictly positive
coefficients fails to be coherent, and the traditional linear pool
fails to be flexibly dispersive.  In view of these limitations, we
investigate parsimonious nonlinear alternatives, including generalized
linear pools (Dawid et al.~1995), the spread-adjusted linear pool,
which has been used successfully in the meteorological literature
(Berrocal, Raftery and Gneiting 2007; Glahn et al.~2009), and the
beta-transformed linear pool (Ranjan and Gneiting 2010), which we
demonstrate to be flexibly dispersive.

Section \ref{sec:examples} turns to a simulation study and data
examples on density forecasts for daily maximum temperature at
Seattle-Tacoma Airport and S\&P 500 returns.  The paper ends in
Section \ref{sec:discussion}, where we summarize our findings from
applied and theoretical perspectives, suggest directions for future
work and, in addition to discussing probabilistic forecasts, hint at
the closely related problem of the fusion of expert judgements that
are expressed in terms of probability distributions.

\section{Combination formulas and aggregation methods}  \label{sec:theory}

In a seminal paper, Murphy and Winkler (1987) proposed a general
framework for the evaluation of point forecasts, which is based on the
joint distribution of the forecast and the observation.  Dawid et
al.~(1995) developed and used a related framework in studying multiple
probability forecasts for a binary event.  Here we respond to the call
of Dawid et al.~(1995, p.~288) for an extension, and we start with an
informal sketch of a fully general approach, in which the observations
take values in just any space.

The most general setting considers the joint distribution of multiple
probabilistic forecasts and the observation on a probability space
$(\Omega, \cA, \myQ)$.  More explicitly, we assume that the elements
of the sample space $\Omega$ can be identified with tupels of the form
\[
(P_1, \ldots, P_k, Y), 
\]
where each of $P_1, \ldots, P_k$ is a probability measure on the
outcome space of the observation, $Y$.  For $i = 1, \ldots, k$, we
require the random probability measure $P_i$ to be measurable with
respect to the sub-$\sigma$-algebra $\cA_i \subseteq \cA$ that encodes
the forecast's information set or information basis, consisting of
data, expertise, theories and assumptions at hand.  The probability
measure $\myQ$ on $(\Omega,\cA)$ specifies the joint distribution of
the probabilistic forecasts and the observation.

In this setting, the probabilistic forecasts $P_1, \ldots, P_k$ might
stem from distinct experts, organizations or statistical models, as
commonly encountered in the practice of forecasting.  In aggregating
them, the theoretically optimal strategy is to combine information
sets, that is, to issue the conditional distribution of the
observation $Y$ given the $\sigma$-algebra $\sigma(\cA_1, \ldots,
\cA_k)$ generated by the information sets $\cA_1, \ldots, \cA_k$.
However, as Dawid et al.~(1995, p.~264) note,
\begin{quote} 
\small
``this ideal will almost always be rendered unattainable, by the
  extent of the data, company confidentiality, or the inability of the
  experts to identify clearly the empirical basis and background
  knowledge leading to their intuitive opinions.``
\end{quote} 
The best that we can hope for in practice is to find the conditional
distribution of the observation $Y$ given the $\sigma$-algebra
$\sigma(P_1, \ldots, P_k)$ generated by the random probability
measures $P_1, \ldots, P_k$.  Of course, it is always true that
\[
\sigma(P_1, \ldots, P_k) \subseteq \sigma(\cA_1, \ldots, \cA_k), 
\]
and in most cases of practical interest the left-hand side constitutes
a substantially lesser information basis than the right-hand side.

\subsection{Prediction spaces}  \label{sec:spaces}

In what follows, we restrict the discussion to the case of a
real-valued observation.  A probabilistic forecast then corresponds to
a Lebesgue-Stieltjes measure on the real line, $\real$, which we
identify with the associated right-continuous cumulative distribution
function (CDF).  We write $\sigma(\cA_1, \ldots, \cA_m)$ and
$\sigma(X_1, \ldots, X_n)$ to denote the $\sigma$-algebra generated by
the families $\cA_1, \ldots, \cA_m$ of subsets of $\Omega$, and the
random variables $X_1, \ldots, X_n$, respectively.  We use the symbol
$\cL$ generically to denote an unconditional or conditional law or
distribution and follow standard conventions in identifying the
sub-$\sigma$-algebras on which we condition.  

We now introduce the key tool of a prediction space, which is a
probability space tailored to the study of combination formulas for
real-valued outcomes, though we allow the case $k = 1$ of a single
probabilistic forecast.

\begin{definition}  \label{def:space}
Let $k \geq 1$ be an integer.  A {\em prediction space}\/ is a
probability space $(\Omega, \cA, \myQ)$ together with
sub-$\sigma$-algebras $\cA_1, \ldots, \cA_k \subseteq \cA$, where the
elements of the sample space $\Omega$ can be identified with tupels
$(F_1, \ldots, F_k, Y, V)$ such that
\begin{itemize}
\item[(P1)] for $i = 1, \ldots, k$, $F_i$ is a CDF-valued random
  quantity that is measurable with respect to the sub-$\sigma$-algebra
  $\cA_i$,\footnote{That is, $\{ F_i(x_j) \in B_j \mbox{ for } j = 1,
    \ldots, n \} \in \cA$ for all finite collections $x_1, \ldots,
    x_n$ of real numbers and $B_1, \ldots, B_n$ of Borel sets.}
\item[(P2)]
  $Y$ is a real-valued random variable, 
\item[(P3)] $V$ is a random variable that is uniformly distributed on
  the unit interval and independent of $\cA_1, \ldots, \cA_k$ and $Y$.
\end{itemize}
\end{definition} 

All subsequent definitions and results are within the prediction space
setting.  Phrases such as {\em almost surely}\/ or {\em with positive
  probability}\/ refer to the probability measure $\myQ$ on $(\Omega,
\cA)$ that determines the joint distribution of the probabilistic
forecasts and the observations.  While (P1) and (P2) formalize the
predictive distributions and the observation, assumption (P3) is
purely technical, allowing us to define a generalized version of the
classical probability integral transform.  The sub-$\sigma$-algebra
$\cA_i$ encodes the information set for the CDF-valued random quantity
$F_i$ which may, but need not, be ideal in the following
sense.\footnote{In a recent comment, Tsyplakov (2011) proposes the
  same terminology.}

\begin{definition}  \label{def:ideal} 
The CDF-valued random quantity $F_i$ is {\em ideal}\/ relative to the
sub-$\sigma$-algebra $\cA_i$ if $F_i = \cL( \, Y \hsp | \hsp\hsp \cA_i
\hsp )$ almost surely.
\end{definition} 

In the subsequent examples we write $\cN(\mu,\sigma^2)$ to denote a
univariate normal distribution with mean $\mu$ and variance
$\sigma^2$, and similarly for the bivariate normal distribution.  We
use the symbols $\Phi$ and $\phi$ to denote the standard normal
cumulative distribution function and density function, respectively. 

\begin{example}[probability forecasts of a binary event]  \label{ex:RG1} 
Here we find a prediction space for the simulation example of Ranjan
and Gneiting (2010), which considers a dichotomous outcome, $Y$.  For
technical consistency later on, we identify a success or occurence
with the outcome $Y = 0$, and a non-success with $Y = 1$.  The
CDF-valued random quantities $F_i$ thus are of the form $F_i(y) = p_i
\hsp \one(y \geq 0) + (1 - p_i) \hsp \one(y \geq 1)$ and can be
identified with the random probability forecasts $p_i$ for a success,
for $i = 1, 2$.  Writing $\omega = (\omega_1, \omega_2, \omega_3,
\omega_4)$ for an elementary event, we define the probability
forecasts
\[
p_1(\omega) = \Phi \! \left( \frac{\omega_1}{\sqrt{1 + \sigma_2^2}} \right) 
\qquad \mbox{and} \qquad
p_2(\omega) = \Phi \! \left( \frac{\omega_2}{\sqrt{1 + \sigma_1^2}} \right) \! , 
\]
where $\sigma_1 > 0$ and $\sigma_2 > 0$ are fixed constants, the
observation $Y(\omega) = \omega_3$, and the auxiliary variable
$V(\omega) = \omega_4$.  To complete the specification of the
prediction space, we let $\Omega = \real \times \real \times \{ 0, 1
\} \times (0,1)$ and let $\cA = \cB_4$ be the corresponding
Borel-$\sigma$-algebra, define $\myQ$ to be the product of
$\cN(0,\sigma_1^2)$, $\cN(0,\sigma_2^2)$ and standard uniform measures
on the first, second and fourth coordinate projections, respectively,
and let
\[
\myQ( B_1 \times B_2 \times \{ 0 \} \times (0,1) ) 
= \frac{1}{\sigma_1^2 \sigma_2^2} \int_{B_1} \int_{B_2} \Phi(\omega_1 + \omega_2) \,
  \phi \! \left( \frac{\omega_1}{\sigma_1} \right) \!   
  \phi \! \left( \frac{\omega_2}{\sigma_2} \right)
  d\lambda(\omega_1) \, d\lambda(\omega_2)    
\]
for Borel sets $B_1, B_2 \subseteq \real$, where $\lambda$ denotes the
Lebesgue measure.  Then $p_1$ is measurable with respect to the
sub-$\sigma$-algebra $\cA_1 = \sigma(\omega_1)$, and $p_2$ is
measurable with respect to $\cA_2 = \sigma(\omega_2)$.  Moreover, by
the arguments in the appendix of Ranjan and Gneiting (2010), $F_1$ is
ideal relative to $\cA_1$, and $F_2$ is ideal relative to $\cA_2$.
\end{example}

Typically, it suffices to consider the joint distribution of the tuple
$(F_1, \ldots, F_k, Y)$, without any need to explicitly specify other
facets of the prediction space.

\begin{table}[t] 
\centering 
\caption{Probabilistic forecasts in Examples \ref{ex:GBR1} and
  \ref{ex:GBR2}.  The observation $Y$ is normal with mean $\mu$ and
  variance 1, where $\mu$ is standard normal.  The random variable
  $\tau$ attains the values $-1$ and $1$ with probability
  $\frac{1}{2}$, independently of $\mu$ and $Y$.  \label{tab:GBR}}

\bigskip
\small
\begin{tabular}{lllll}
\hline
\hline
Forecast & Predictive Distribution \rule{0mm}{4.5mm} \\
\hline
Perfect        & $F_1 = \cN(\mu,1)$ \rule{0mm}{4.5mm} \\
Climatological & $F_2 = \cN(0,2)$ \rule{0mm}{4.5mm} \\
Unfocused      & $F_3 = \frac{1}{2} \left( \hsp \cN(\mu,1) + \cN(\mu+\tau,1) \right)$ \rule{0mm}{4.5mm} \\
Sign-reversed  & $F_4 = \cN(-\mu,1)$ \rule{0mm}{4.5mm} \\
\hline
\end{tabular}
\end{table} 

\begin{example}[density forecasts]  \label{ex:GBR1} 
To define a prediction space, let
\[
Y \mid \mu \sim \cN(\mu,1) \qquad \mbox{where} \qquad \mu \sim \cN(0,1),  
\]
and let $\tau$ attain the values 1 and $-1$ with equal probability,
independently of $\mu$ and $Y$.  Table \ref{tab:GBR} places the
density forecasts in the simulation example of Gneiting, Balabdaoui
and Raftery (2007) in this setting.  The perfect forecast is ideal
relative to the sub-$\sigma$-algebra generated by $\mu$.  The
climatological forecast is ideal relative to the trivial
sub-$\sigma$-algebra.  The unfocused and sign-biased forecasts are not
ideal, as we will see in Example \ref{ex:GBR2} below.
\end{example}

\subsection{Combination formulas}  \label{sec:formulas}

As noted, in aggregating predictive cumulative distribution functions,
the theoretically optimal strategy is to combine information sets,
that is, to issue the conditional distribution of the observation $Y$
given the $\sigma$-algebra $\sigma(\cA_1, \ldots, \cA_k)$ generated by
the information sets $\cA_1, \ldots, \cA_k$.  However, information
aggregation often is not feasible in practice, when individual sources
of expertise reveal predictive distributions, rather than information
sets.  What we can realistically aim at is to model the conditional
distribution of the observation $Y$ given the $\sigma$-algebra
generated by the predictive cumulative distribution functions, 
namely 
\[
G = \cL \, ( \, Y \hsp\hsp | \hsp\hsp F_1, \ldots, F_k \hsp ), 
\]
where we define 
\[
\cL \, ( \, Y \hsp\hsp | \hsp\hsp F_1, \ldots, F_k \hsp )
= \cL \, ( \, Y \hsp\hsp | \hsp\hsp F_i(x) : i = 1, \ldots, k, \, x \in \qq ), 
\] 
with $\qq$ being the set of the rational numbers. 

\begin{table}[t] 
\centering 
\caption{Some classes of fixed, non-random cumulative distribution
  functions, where the subscript refers to an interval $\myI \subseteq
  \real$.  In the case of Bernoulli measures, we identify a success
  with 0 and a non-success with 1, so that the corresponding
  cumulative distribution function has jump discontinuities at these
  values, and otherwise is constant.
\label{tab:F}}

\bigskip
\small
\begin{tabular}{ll}
\hline
\hline
Class & Characterization of the Members \rule{0mm}{4.25mm} \\
\hline
$\cF_\myI$   & support in $\myI$ \rule{0mm}{4.25mm} \\
$\cF_\myI^+$ & support in $\myI$; strictly increasing on $\myI$ \\
$\cC_\myI$   & support in $\myI$; continuous \rule{0mm}{4mm} \\
$\cC_\myI^+$ & support in $\myI$; continuous; strictly increasing on $\myI$ \\
$\cD_\myI$   & support in $\myI$; admits Lebesgue density \rule{0mm}{4mm} \\
$\cD_\myI^+$ & support in $\myI$; admits Lebesgue density; strictly increasing on $\myI$ \\
$\cB$       & Bernoulli measure \rule{0mm}{4mm} \\
$\cB^+$     & Bernoulli measure with nondegenerate success probability \\
\hline 
\end{tabular}
\end{table} 

In practice one resorts to empirical combination formulas.
Specifically, let $\cF$ be a class of fixed, non-random cumulative
distribution functions such that $F_1, \ldots, F_k \in \cF$ almost
surely.  For example, if we are concerned with density forecasts on on
the real line $\real$, we consider the class $\cD_\real$ of the
cumulative distribution functions that admit a Lebesgue density.
Further classes $\cF$ of interest are listed in Table \ref{tab:F}.
A {\em combination formula}\/ then is a mapping of the form
\begin{equation}  \label{eq:G} 
G : \cF^k = \underbrace{\: \cF \times \cdots \times \cF}_{k \; \mbox{\footnotesize times}} \to \cF, 
\qquad (F_1, \ldots, F_k) \mapsto G(F_1, \ldots, F_k).  
\end{equation}
We allow the case $k = 1$, where the mapping may
provide calibration and dispersion adjustments for a single predictive
distribution, as described later on in Sections \ref{sec:BLP} and
\ref{sec:discussion}.

Here our interest is in the case $k \geq 2$, where Dawid et al.~(1995)
introduced the notion of a coherent combination formula in the context
of probability forecasts of a binary event.  We now define a more
general notion that applies to probabilistic forecasts of general,
real-valued quantities.

\begin{definition}  \label{def:coherent} 
Suppose that $k \geq 2$.  The combination formula $G : \cF^k \to \cF$
is {\em coherent}~relative to the class $\cF$ if there exists a
prediction space $(\Omega, \cA, \myQ)$ with sub-$\sigma$-algebras
$\cA_1, \ldots, \cA_k \subseteq \cA$ and CDF-valued random quantities
$F_1, \ldots, F_k$ such that
\begin{itemize}
\item[(C1)] for $i = 1, \ldots, k$, $F_i = \cL(Y \hsp\hsp | \hsp\hsp
  \cA_i) \in \cF$ almost surely,
\item[(C2)] for $i \not= j$, $F_i \not= F_j$ with positive
  probability,
\item[(C3)] $\cL( \hsp Y \hsp\hsp | \hsp\hsp F_1, \ldots F_k ) =
  G( F_1, \ldots, F_k)$ almost surely.
\end{itemize} 
\end{definition} 

It is important to note that condition (C1) has two requirements, the
first being that $F_i = \cL(Y \hsp\hsp | \hsp\hsp \cA_i)$ is ideal,
and the second that $F_i \in \cF$ almost surely.  Condition (C2)
excludes trivial cases, while (C3) requests the aggregated cumulative
distribution function to be ideal relative to the $\sigma$-algebra
generated by the full set $F_1, \ldots, F_k$ of the components.  Note
that the smaller the class $\cF$, the stronger a statement about
coherence.  Conversely, the larger the class $\cF$, the stronger a
statement about lack of coherence.  

Conceptually, a coherent combination formula is compatible with fully
informed and theoretically literate decision making, whereas an
incoherent combination formula is not.  In this light, coherence is an
appealing property of a combination formula.  Nevertheless, the
applied relevance of the notion of coherence is limited, as
probabilistic forecasts are hardly ever ideal in practice, which may
invite, or even necessitate, the use of flexible classes of
potentially incoherent combination formulas.  Motivated by these
applied perspectives, we now turn to a discussion of the important
notions of calibration and dispersion.

\subsection{Calibration and dispersion}  \label{sec:dispersion}

If $F$ denotes a fixed, non-random predictive cumulative distribution
function for an observation $Y$, the probability integral transform is
the random variable $Z_F = F(Y)$.  It is well known that if $F$ is
continuous and $Y \sim F$ then $Z_F$ is standard uniform (Rosenblatt
1952).  If the more general, randomized version of the probability
integral transform proposed by Brockwell (2007) is used, the
uniformity result applies to arbitrary, not necessarily continuous,
but still fixed, non-random cumulative distribution functions.  In the
prediction space setting, we need the following, further extension
that allows for $F$ to be a CDF-valued random quantity.

\smallskip
\begin{definition}  
In the prediction space setting, the random variable
\[
Z_F = {\textstyle \lim_{\, y \uparrow Y}} F(y) 
      + V \left( \rule{0mm}{3.5mm} F(Y) - {\textstyle \lim_{\, y \uparrow Y}} F(y) \right) 
\]
is the {\em probability integral transform}\/ of the CDF-valued random
quantity $F$.
\end{definition}

In a nutshell, the probability integral transform is the value that
the predictive cumulative distribution function attains at the
observation, with suitable adaptions at points of discontinuity.
The probability integral transform takes values in the unit interval,
and so the possible values of its variance are constrained to the
closed interval $[0,\frac{1}{4}]$.  A variance of $\frac{1}{12}$
corresponds to a uniform distribution and continues to be the most
desirable, as evidenced by Theorem \ref{th:ideal} below.

We are now ready to define and study notions of calibration and
dispersion.  In doing so, we use the terms CDF-valued random quantity
and forecast interchangeably.

\smallskip
\begin{definition}  \label{def:calibration} 
In the prediction space setting, let $F$ and $G$ be CDF-valued random
quantities with probability integral transforms $Z_F$ and $Z_G$.
\begin{enumerate}
\item[\rm (a)] The forecast $F$ ist {\em marginally calibrated}\/ if
  $\myE_{\, \myQ} [F(y)] = \myQ(Y \leq y)$ for all $y \in \real$.
\item[\rm (b)] The forecast $F$ is {\em probabilistically
  calibrated}\/ if its probability integral transform $Z_F$ is
  uniformly distributed on the unit interval.
\item[\rm (c)] The forecast $F$ is {\em overdispersed}\/ if $\var(Z_F)
  < \frac{1}{12}$, {\em neutrally dispersed}\/ if $\var(Z_F) =
  \frac{1}{12}$, and {\em underdispersed}\/ if $\var(Z_F) >
  \frac{1}{12}$.
\item[\rm (d)] The forecast $F$ is {\em at least as
  dispersed}\/ as the forecast $G$ if $\var(Z_F) \leq \var(Z_G)$.  It
  is {\em more dispersed}\/ than $G$ if $\var(Z_F) < \var(Z_G)$.
\item[\rm (e)] The forecast $F$ is {\em regular}\/ if the distribution
  of $Z_F$ is supported on the unit interval.
\end{enumerate}
\end{definition}

In the defining equality $\myE_{\, \myQ} [F(y)] = \myQ(Y \leq y)$ for
marginal calibration, the left-hand side depends on the law of the
predictive distribution, whereas the right-hand side depends on the
law of the observation.  In parts (c) and (d) of Definition
\ref{def:calibration}, we define dispersion in terms of the variance
of the probability integral transform, thus involving the joint law of
the predictive distribution and the observation.  In contrast, the
spread of the predictive distribution itself is a measure of sharpness
that does not consider the observation.

Our current setting of prediction spaces differs from, but relates
closely to, the approach of Gneiting et al.~(2007), who studied
notions of calibration from a prequential perspective.  Specifically,
if the CDF-valued random quantity $F$ is probabilistically calibrated
in the sense of Definition \ref{def:calibration} and we sample from
the joint law of $F$ and $Y$, the resulting sequence is
probabilistically calibrated in the sense of Gneiting et al.~(2007).
An analogous statement applies to marginal calibration.

Returning to the prediction space setting, the following result is
immediate.

\begin{proposition}  \label{prop:neutrally}
A probabilistically calibrated forecast is neutrally dispersed
and regular.
\end{proposition}

The converse is not necessarily true, in that a forecast which is
neutrally dispersed need not be calibrated nor regular.  However,
an ideal forecast is always calibrated. 

\begin{theorem}  \label{th:ideal}
A forecast that is ideal relative to a $\sigma$-algebra ist both
marginally calibrated and probabilistically calibrated.
\end{theorem}

{\em Proof.} Suppose that $F = \cL(Y \hsp | \hsp \cA_0)$ is ideal
relative to the $\sigma$-algebra $\cA_0$, so that $F(y) = \myQ( Y \leq
y \hsp\hsp | \hsp\hsp \cA_0)$ almost surely for all $y \in \real$.
Then
\[
\myE_\myQ [F(y)]
= \myE_\myQ \! \left[ \hsp \myQ( Y \leq y \hsp\hsp | \hsp\hsp \cA_0) \right]
= \myE_\myQ \, \myE_\myQ \left[ \hsp\hsp \one \hsp (Y \leq y) \hsp\hsp | 
  \hsp\hsp \cA_0 \hsp \right] 
= \myQ(Y \leq y),
\]
where $\one$ denotes an indicator function, thereby proving the
statement about marginal calibration.  Turning to probabilistic
calibration, let $\myQ_0$ denote the marginal law of $Y$ under $\myQ$,
so that $Z_F = \myQ_0((-\infty,Y) \hsp\hsp | \hsp\hsp \cA_0) + V \,
\myQ_0(\{Y\} \hsp\hsp | \hsp\hsp \cA_0)$ and
\[
\myQ(Z_F \leq z) 
= \myE_\myQ \, \myE_\myQ \left[ \hsp\hsp \one \hsp (Z_F \leq z) \hsp\hsp | 
  \hsp\hsp \cA_0 \hsp \right] 
= z
\]
for $z \in (0,1)$, where the final equality uses the key result of
Brockwell (2007).  \done

\begin{example} \label{ex:GBR2} 
We revisit the density forecasts in Example \ref{ex:GBR1} and Table
\ref{tab:GBR}.  The perfect forecast and the climatological forecast
are ideal, and so by Theorem \ref{th:ideal} they are both
probabilistically calibrated and marginally calibrated.  Arguments
nearly identical to those in Gneiting et al.~(2007) show that the
unfocused forecaster is probabilistically calibrated but not
marginally calibrated, and that the sign-biased forecaster is
marginally calibrated but not probabilistically calibrated.  Hence,
there is no sub-$\sigma$-algebra or information set relative to which
the unfocused or the sign-biased forecaster is ideal.
\end{example}

Dawid (1984), Diebold et al.~(1998), Gneiting et al.~(2007) and Czado
et al.~(2009), among others, have argued powerfully that probabilistic
calibration is a critical requirement for probabilistic forecasts that
take the form of predictive cumulative distribution functions, with
Theorem \ref{th:ideal} lending further support to this approach.
Indeed, checks for the uniformity of the probability integral
transform have formed a cornerstone of density forecast evaluation.
In practice, one observes a sample $\{ (F_{1j}, \ldots, F_{kj}, y_j) :
j = 1, \ldots, J \}$ from the joint distribution of the probabilistic
forecasts and the observation, and the uniformity of the probability
integral transform is assessed empirically.  The prevalent way of
doing this is by plotting histograms of the probability integral
transform values for the various forecasting methods, which show the
corresponding frequency distribution over an evaluation or test set.
U-shaped histograms correspond to underdispersed predictive
distributions with prediction intervals that are too narrow on
average, while hump or inverse U-shaped histograms indicate
overdispersed predictive distributions.  Formal tests of uniformity
can also be employed; for a review, see Corradi and Swanson (2006).

\begin{figure}[p]
\centering
\includegraphics[width = 0.5\textwidth]{\mypath 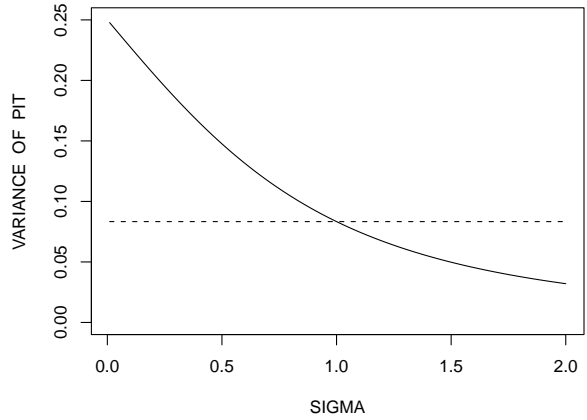}
\caption{The variance (\ref{eq:var}) of the probability integral
  transform $Z_\sigma = F_\sigma(Y)$ for the predictive distribution
  $F_\sigma$ in Example \ref{ex:dispersion} as a function of the
  predictive standard deviation, $\sigma$.  The dashed horizontal line
  at $\frac{1}{12}$ indicates a neutrally dispersed
  forecast.  \label{fig:variance}}
\end{figure}

\begin{figure}[p]
\centering
\includegraphics[width = \textwidth]{\mypath 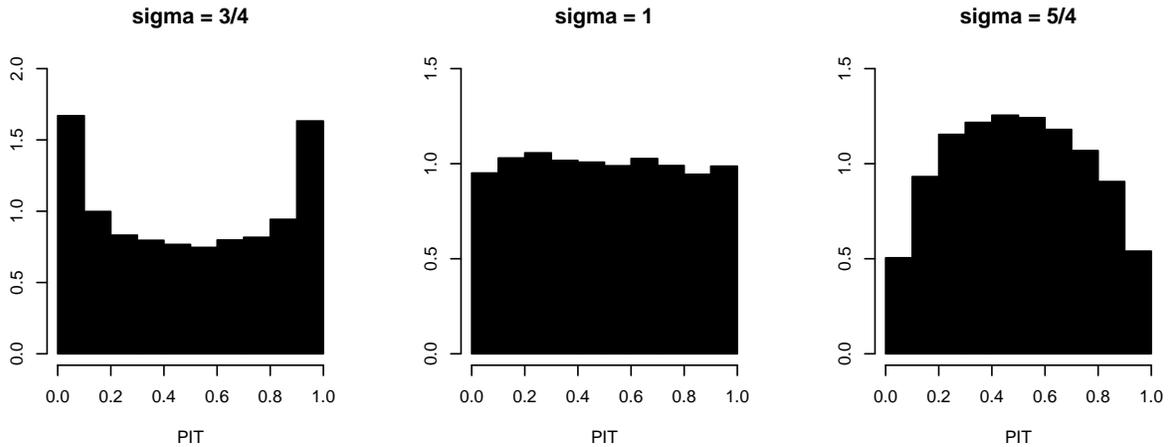}
\caption{Probability integral transform histograms for the predictive
  distribution $F_\sigma$ in Example \ref{ex:dispersion}, where
  $\sigma = \frac{3}{4}$ (underdispersed), $\sigma = 1$ (neutrally
  dispersed and calibrated) and $\sigma = \frac{5}{4}$
  (overdispersed).
  \label{fig:PIT}}
\end{figure}

\begin{example}  \label{ex:dispersion}
Let $Y = X + \epsilon$, where $X$ and $\epsilon$ are independent,
standard normal random variables, and consider the Gaussian predictive
distribution $F_\sigma = \cN(X,\sigma^2)$.  A stochastic domination
argument, the details of which we give in Appendix A, shows that
$F_\sigma$ is underdispersed if $\sigma < 1$, neutrally dispersed if
$\sigma = 1$ and overdispersed if $\sigma > 1$.  If $\sigma = 1$ then
$F_\sigma$ is ideal and thus both marginally calibrated and
probabilistically calibrated.  A more detailed calculation, which is
also given in Appendix A, shows that the probability integral
transform $Z_\sigma = F_\sigma(Y)$ satisfies
\begin{equation}  \label{eq:var}
\var(Z_\sigma)
= 2 \int_0^1 z \left( 1 - \Phi( \sigma \hsp (\Phi^{-1}(z))) \right) \dd z
  - \left( \, \int_0^1 \left( 1 - \Phi( \sigma \hsp (\Phi^{-1}(z))) \right) \dd z \right)^{\! 2}.
\end{equation}
In Figure \ref{fig:variance} we plot $\var(Z_\sigma)$ as a function of
the predictive standard deviation, $\sigma$.  Figure \ref{fig:PIT}
shows probability integral transform histograms for a Monte Carlo
sample of size $10,000$ from the joint distribution of the observation
$Y$ and the forecasts $F_\sigma$, where $\sigma = \frac{3}{4}$, $1$
and $\frac{5}{4}$.  The histograms are U-shaped, uniform, and inverse
U-shaped, reflecting underdispersion, neutral dispersion and
calibration, and overdispersion, respectively.
\end{example} 

Recall from Example \ref{ex:RG1} that in the case of a binary outcome
$Y$, we identify a CDF-valued random quantity $F(y) = p \hsp \one(y
\geq 0) + (1 - p) \hsp \one(y \geq 1)$ with the probability forecast
$p$ for a success, that is, $Y = 0$.  The extant literature, including
Schervish (1989) and Ranjan and Gneiting (2010) and the references
therein, calls $p$ calibrated if
\begin{equation}  \label{eq:cond.cal} 
\myQ( \hsp Y = 0 \hsp\hsp | \, p \hsp ) = p \qquad \mbox{almost surely}. 
\end{equation} 
Here we refer to this property as {\em conditional calibration}.
Perhaps surprisingly, our next result shows that if the outcome is
binary, the notions of probabilistic calibration and conditional
calibration are equivalent.  For an illustrating example, see Appendix
C.

\begin{theorem}  \label{th:binary} 
Consider a prediction space\/ $(\Omega, \cA, \myQ)$ with a binary
outcome\/ $Y$, where\/ $Y = 0$ corresponds to a success and\/ $Y = 0$
to a failure, and a CDF-valued random quantity\/ $F(y) = p \, \one(y
\geq 0) + (1 - p) \hsp \one(y \geq 1)$, which can be identified with
the probability forecast\/ $p$ for a success.  Then the following
statements are equivalent:
\begin{enumerate}
\item[\rm (i)] The forecast\/ $F$ is probabilistically calibrated,
  that is, its probability integral transform\/ $Z_F$ is uniformly
  distributed on the unit interval.
\item[\rm (ii)] The probability forecast\/ $p$ is conditionally
  calibrated, that is, $\myQ( \hsp Y = 0 \hsp\hsp | \, p \hsp ) = p$
  almost surely.
\item[\rm (iii)] The forecast\/ $F$ is ideal relative to the
  $\sigma$-algebra generated by the probability forecast\/ $p$. 
\end{enumerate}
\end{theorem} 

{\em Proof.}  It is clear that (ii) and (iii) are equivalent, and by
Theorem \ref{th:ideal} the statement (iii) implies (i).  To conclude
the proof, we show that statement (i) implies (ii).  To this end,
suppose that the forecast $F$ is probabilistically calibrated.  By
standard properties of conditional expectations, there exists a
measurable function $q : [0,1] \to [0,1]$ such that $\myQ( \hsp Y = 0
\hsp\hsp | \, p \hsp ) = q( \hsp p)$ almost surely.  Let $H$ denote
the marginal law of $p$ under $\myQ$.  If $H$ has a point mass at 0 or
1, it is readily seen that $q(0) = 0$ or $q(1) = 1$, respectively.

A version of the conditional density $u(z \hsp | \hsp x)$ of the
probability integral transform $Z_F$ given that $p = x \in [0,1]$
satisfies $u(z \hsp | \hsp x) = (1-q(x))/(1-x)$ for $z \in [0,1-x]$
and $u(z \hsp | \hsp x) = q(x)/x$ for $z \in (1-x,1]$.  The marginal
  density $u$ of $Z_F$ is standard uniform, so that
\[
u(z + \delta) - u(z) 
= \int_{[0,1]} \!
  \left( \rule{0mm}{3.5mm} \hsp u(z + \delta \hsp | \hsp x) - u(z \hsp | \hsp x) \right) \dd H(x) 
= \int_{(1-z-\delta,1-z]} \frac{q(x)-x}{x \hsp (1-x)} \, \dd H(x) = 0 
\]
for Lebesgue-almost all $z \in (0,1)$ and $\delta \in (0, 1-y)$.  
Let $0 < a < b < 1$, and consider the signed measure defined by 
\[
\mu(A) = \int_A \, \frac{q(x)-x}{x \hsp (1-x)} \, \dd H(x)
\]
for Borel sets $A \subseteq [a,b]$.  We have just shown that $\mu(A) =
0$ for all intervals $(c,d] \subseteq [a,b]$, except possibly for $c$
  or $d$ in a Lebesgue null set.  Since the family of such intervals
  generates the Borel-$\sigma$-algebra, this is only possible if $\mu$
  is the null measure, so that $\mu(A) = 0$ for all Borel sets $A$ in
  $[a,b]$.  In particular, $\mu(A) = 0$ for $A = \{ x \in [a,b] : q(x)
  > x \}$ and $A = \{ x \in [a,b] : q(x) < x \}$, which implies that
  $q(x) = x$ almost surely with respect to the restriction of $H$ to
  $[a,b]$.  To summarize, we have shown that $q(x) = x$ almost surely
  with respect to $H$, whence $\myQ( \hsp Y = 0 \hsp\hsp | \, p \hsp )
  = p$ almost surely with respect to $\myQ$, as desired.  \done

\begin{table}[t]
\centering
\caption{Example of a probabilistically calibrated, but not
  auto-calibrated CDF-valued random quantity $F$ for a ternary outcome
  $Y$.
\label{tab:ternary}} 
\footnotesize
\medskip
\begin{tabular}{cccccccc}
\hline
\hline
$\myQ$-probability \rule{0mm}{4mm} 
& \multicolumn{4}{c}{$F(x)$} & \multicolumn{3}{c}{$\myQ(Y = i \hsp\hsp | \hsp F)$} \\ 
\rule{0mm}{4mm} & $x < 0$ & $0 \leq x < 1$ & $1 \leq x < 2$ & $x \geq 2$ & 
$i = 0$ & $i = 1$ & $i = 2$ \\  
\hline
$\frac{1}{2}$ & 0 & $\frac{1}{2}$ & 1 & 1 & $\frac{3}{4}$ & $\frac{1}{4}$ & 0 \rule{0mm}{4mm} \\
$\frac{1}{2}$ & 0 & $\frac{1}{2}$ & $\frac{3}{4}$ & 1 & $\frac{1}{4}$ & $\frac{3}{8}$ & $\frac{3}{8}$ 
\rule{0mm}{4mm} \vphantom{$\frac{1}{2_f}$} \\
\hline
\end{tabular}
\smallskip
\end{table}

\bigskip
Theorem \ref{th:binary} draws a connection from the probability
integral transform histogram to the {\em reliability diagram}\/ or
{\em calibration curve}, which is the key diagnostic tool for
assessing the calibration of probability forecasts for a binary event
(Dawid 1986; Murphy and Winkler 1992; Ranjan and Gneiting 2010).  A
reliability diagram plots conditional event frequencies against binned
forecast probabilities, with deviations from the diagonal indicating
violations of the conditional calibration condition
(\ref{eq:cond.cal}).  For non-binary outcomes $Y$ and the natural
generalization of the conditional calibration criterion, namely the
{\em auto-calibration}\/ property
\[
\cL( \hsp Y \hsp\hsp | \, F \hsp ) = F \qquad \myQ\mbox{-almost surely}  
\]
introduced by Tsyplakov (2011) in a comment on Mitchell and Wallis
(2010), the equivalence to probabilistic calibration fails, as
demonstrated in Table \ref{tab:ternary} for a ternary outcome.  In
general, auto-calibration is a much stronger requirement than
probabilistic calibration, and there is a small but colorful, growing
strand of recent literature that addresses its applied aspects for
continuous outcomes (Hamill 2001; Mason et al.~2007; Held, Rufibach
and Balabdaoui 2010; Br\"ocker, Siegert and Kantz 2011; Mason et
al.~2011).

\subsection{Aggregation methods}  \label{sec:aggregation}

An {\em aggregation method}\/ is a family 
\[
\cG = \{ G_\theta : \theta \in \Theta \}
\]
of combination formulas of the form 
\begin{equation}  \label{eq:G.theta} 
G_\theta : \cF^k = \cF \times \cdots \times \cF \to \cF, 
\qquad (F_1, \ldots, F_k) \mapsto G_\theta(F_1, \ldots, F_k).  
\end{equation} 
that share a common value of $k$ and a common class $\cF$ of fixed,
non-random cumulative distribution functions.  For example, if $\cG$
is the traditional linear pool, we can take $\cF$ to be any convex
class of cumulative distribution functions, and we may identify the
index set $\Theta$ with the unit simplex in $\real^k$.

The extant literature studies individual combination formulas by
examining whether or not they possess certain analytic
characteristics, such as the strong setwise function and external
Bayes properties (McConway 1981; Genest 1984a, 1984b; Genest and Zidek
1986; Genest, McConway and Schervish 1986; French and R\'{i}os Insua
2000).  In contrast, we share the recent perspective of Hora (2010)
and put the focus on calibration and dispersion.  In particular, we
study aggregation methods in terms of the behavior of the probability
integral transform under the corresponding family $\cG = \{ G_\theta :
\theta \in \Theta \}$ of combination formulas.  As noted, the
probability integral takes values in the unit interval, and so the
possible values of its variance lie between 0 and $\frac{1}{4}$.  The
value $\frac{1}{12}$ corresponds to neutral dispersion and is the most
desirable.

Following French and R\'{i}os Insua (2000, p.~113) we say that the
combination formula $G_\theta$ is {\em anonymous}\/ if
\[
G_\theta(F_{\pi(1)}, \ldots, F_{\pi(k)}) = G_\theta(F_1, \ldots, F_k) 
\]
for all $F_1, \ldots, F_k \in \cF$ and all permutations $\pi$ on $k$
elements.  With this, we are ready to define notions of flexibility for
aggregation methods.

\begin{definition}  \label{def:dispersive}
Consider a family $\cG = \{ G_\theta : \theta \in \Theta \}$ of
combination formulas of the form (\ref{eq:G.theta}) that share a
common $k \geq 2$ and a common class $\cF$ of fixed, non-random
cumulative distribution functions.  
\begin{itemize}
\item[\rm (a)] The aggregation method $\cG$ is {\em flexibly
  dispersive}\/ relative to the class $\cF$ if for all $F_0 \in \cF$
  and $F_1, \ldots, F_k \in \cF$ there exists a parameter value
  $\theta \in \Theta$ such that if $\cL(Y) = F_0$ then
  $G_\theta(F_1,\ldots,F_k)$ is neutrally dispersed.
\item[\rm (b)] The aggregation method $\cG$ is {\em exchangeably
  flexibly dispersive}\/ relative to the class $\cF$ if for all $F_0
  \in \cF$ and $F_1, \ldots, F_k \in \cF$ there exists a parameter
  value $\theta \in \Theta$ such that $G_\theta$ is anonymous and if
  $\cL(Y) = F_0$ then $G_\theta(F_1,\ldots,F_k)$ is neutrally
  dispersed.
\end{itemize}
\end{definition}

The applied relevance of the definitions is appreciated as follows.
Suppose that the aggregation method $\cG$ is flexibly dispersive
relative to $\cF$.  Then, given any marginal law $F_0 \in \cF$ for the
observation $Y$ and any collection $F_1, \ldots, F_k$ of probabilistic
forecasts for $Y$, we can find a combination formula $G_\theta \in
\cG$ such that the aggregated predictive distribution, namely
$G_\theta(F_1, \ldots, F_k)$, is neutrally dispersed.  If $\cG$ is
exchangeably flexibly dispersive, we can do so while treating $F_1,
\ldots, F_k$ exchangeably, which is a frequent requirement in the
practice of the combination of expert judgements (Jouini and Clemen
1986).

Note that a positive statement about flexible dispersivity is the
stronger, the larger the class $\cF$.  Conversely, a statement about
the lack of flexible dispersivity is the stronger, the smaller the
class $\cF$.

\section{Linear and nonlinear aggregation methods}  \label{sec:methods}

In this section we investigate specific combination formulas and
aggregation methods from the perspectives of coherence, calibration
and dispersion.  First, we consider the traditional linear pool, then
we move on to discuss non-linear ramifications, namely generalized
linear pools, the spread-adjusted linear pool and the beta transformed
linear pool.  Note that the linearity of a combination formula for
cumulative distribution functions is preserved in the corresponding
formula for aggregating densities or event probabilities.  In
contrast, the functional form of a nonlinear combination formula
changes if it is expressed in terms of densities or event
probabilities.

\subsection{Linear pool}  \label{sec:linear}

We proceed to state and prove a simple but powerful result about
linear combination formulas that generalizes earlier findings by Hora
(2004) and Ranjan and Gneiting (2010).  The gist of the statement is
that dispersion tends to increase under linear aggregation.

\smallskip
\begin{theorem}  \label{th:linear}
In the prediction space setting, suppose that\/ $k \geq 2$ and
consider any linearly combined probabilistic forecast\/ $F =
\sum_{i=1}^k w_i \hsp F_i$ with weights\/ $w_1, \ldots, w_k$ that are
strictly positive and sum to\/ $1$.  For\/ $i \not= j$, suppose that\/
$F_i \not= F_j$ with positive probability.  Then the following holds:
\begin{itemize}
\item[\rm (a)] The linearly combined forecast\/ $F$ is at least as
  dispersed as the least dispersed of the components\/ $F_1, \ldots,
  F_k$.
\item[\rm (b)] If the component forecasts\/ $F_1, \ldots, F_k$ are
  regular, then $F$ is more dispersed than the least dispersed of the
  components.
\item[\rm (c)] If the components are neutrally dispersed and regular,
  then $F$ is overdispersed.
\end{itemize}
\end{theorem}

{\em Proof.} \ For $i = 1, \ldots, k$, let $Z_i$ denote the
probability integral transform of $F_i$.  The probability integral
transform of $F = \sum_{i=1}^k w_i \hsp F_i$ is $Z = \sum_{i=1}^k w_i
\hsp Z_i$, whence
\[
\var(Z) = \sum_{i=1}^k \sum_{j=1}^k w_i w_j \, \cov(Z_i,Z_j)
        \leq \sum_{i=1}^k w_i \sum_{j=1}^k w_j \! \left( \max_{1 \leq i \leq k} \var(Z_i) \! \right)
        = \max_{1 \leq i \leq k} \var(Z_i),
\]
which demonstrates part (a).  To prove part (b) suppose, for a
contradiction, that $F_1, \ldots, F_k$ are regular and $\var(Z) =
\max_{1 \leq i \leq k} \var(Z_i)$.  Then $Z_i$ and $Z_j$ are perfectly
correlated for $i, j = 1, \ldots, k$, and we conclude that there exist
constants $a_{ij} > 0$ and $b_{ij} \in \real$ such that $Z_i = a_{ij}
Z_j + b_{ij}$ almost surely.  By the assumption of regularity, $Z_i$
and $Z_j$ are supported on the unit interval, whence $a_{ij} = 1$ and
$b_{ij} = 0$.  Therefore, $Z_i = Z_j$ almost surely, contrary to the
assumption that $F_i \not= F_j$ with positive probability.  Part (c)
concerns the special case of part (b) in which $\var(Z_i) =
\frac{1}{12}$ for $i = 1, \ldots, k$.  \done

\bigskip
As noted, Theorem \ref{th:linear} yields various extant results as
corollaries.  For instance, Hora (2004) applied Fourier analytic tools
to show that if two distinct density forecasts are probabilistically
calibrated, then any nontrivial linear combination is uncalibrated,
which is an immediate consequence of part (c).  However, the statement
of part (c) is considerably stronger, in that it substitutes the
weaker condition of neutral dispersion and regularity for the
assumption of probabilistic calibration, allows for any number $k \geq
2$ of components, allows for cumulative distribution functions rather
than the special case of densities, and exposes the direction of the
deviation, in that the linearly combined forecast is overdispersed.
Each of the four facets is useful in practice.  For instance, there
are real data situations where density forecasts are approximately
neutrally dispersed and regular, but clearly not calibrated.  Discrete
and mixed discrete-continuous predictive cumulative distribution
functions also occur frequently in practice, such as in quantitative
precipitation forecasting (Sloughter et al.~2007) and for count data
(Czado et al.~2009).  Finally, the tendency to increase dispersion
helps explain the success of linear pooling in applications, where the
component distributions are frequently underdispersed.  For a
prominent example, see Table 10 of Hoeting et al.~(1999).

The next theorem nests extant results in the case of probability
forecasts for a binary event, which were proved by Dawid et al.~(1995)
and Ranjan and Gneiting (2010).

\begin{theorem}  \label{th:linear.incoherent} 
If\/ $k \geq 2$, any linear combination formula with strictly positive
weights fails to be coherent relative to the class\/ $\cF_\real$.
\end{theorem}

{\em Proof.}  Suppose, for a contradiction, that the combination
formula $G(F_1, \ldots, F_k) = F = \sum_{i=1}^k w_i F_i$ with strictly
positive weights that sum to 1 is coherent relative to $\cF_\real$.
Then there exists a prediction space $(\Omega, \cA, \myQ)$ with
sub-$\sigma$-algebras $\cA_1, \ldots, \cA_k \subseteq \cA$ and
CDF-valued random quantities $F_1, \ldots, F_k$ such that properties
(C1), (C2) and (C3) of Definition \ref{def:coherent} hold.  By
property (C1) and Theorem \ref{th:ideal}, the components $F_1, \ldots,
F_k$ are probabilistically calibrated, and by Proposition
\ref{prop:neutrally} they are neutrally dispersed and regular.  Thus,
part (c) of Theorem \ref{th:linear} applies and the linearly
aggregated forecast $F$ is overdispersed.  However, by property (C3),
Theorem \ref{th:ideal} and Proposition \ref{prop:neutrally} $F$ is
neutrally dispersed, for the desired contradiction.  \done

\bigskip
Thus far, we have considered individual linear combination formulas.
The following result views the traditional linear pool as an
aggregation method $\cG = \{ G_\theta : \theta \in \Theta \}$, where
we may identfy the parameter space $\Theta = \Delta_{k-1}$ with the
unit simplex in $\real^k$.  We state the theorem relative to the full
class $\cF_\real$, even though it remains valid relative to much
smaller classes.

\begin{theorem}
If\/ $k \geq 2$, the linear pool fails to be flexibly dispersive
relative to the class\/ $\cF_\real$.
\end{theorem}

{\em Proof.}  In view of part (a) of Theorem \ref{th:linear} it
suffices to find an $F_0 \in \cF_\real$ and distinct $F_1, \ldots, F_k
\in \cF_\real$, each of which is an overdispersed as a probabilistic
forecast for an observation $Y$ with $\cL(Y) = F_0$.  For example, we
can take $F_0$ to be standard normal and $F_i$ to be normal with mean
zero and variance $i + 1$ for $i = 1, \ldots, k$.  \done

\subsection{Generalized linear pools}  \label{sec:GLP}

Dawid et al.~(1995) introduced and studied {\em generalized linear}\/
combination formulas for combining probability forecasts of a binary
event.  Here we apply the approach to cumulative distribution
functions, where we obtain combination formulas of the form
\begin{equation}  \label{eq:GLP} 
G(y) = h^{-1} \! \left( \sum_{i=1}^k w_i \, h(F_i(y)) \right) 
\qquad \mbox{or} \qquad 
h(G(y)) = \sum_{i=1}^k w_i \, h(F_i(y)), 
\end{equation} 
where $h$ is a continuous and strictly monotone link function.  It is
interesting to note the formal resemblance to Archimedean copulas
(Genest and MacKay 1986; McNeil and Ne\u{s}lehov\'a 2009), as observed
and investigated by Jouini and Clemen (1996).

\begin{table}[t] 
\centering 
\caption{Specifics of the generalized linear combination formula in
  equation (\ref{eq:GLP}).  The table states assumptions on the
  weights, $w_1, \ldots, w_k$, and instances of classes $\cF$, such
  that the combination formula maps $\cF^k$ into $\cF$.  The
  conditions depend on the domain and the range of the link function,
  $h$, which we assume to be continuous and strictly
  monotone. \label{tab:h}}

\bigskip
\small
\begin{tabular}{llllll}
\hline
\hline
Type & Domain & Range & Weights & Class $\cF$ & Example \rule{0mm}{4.5mm} \\
\hline
A & $[0,1]$ & any           & $w_i \geq 0$; $\sum_{i=1}^k w_i = 1$ & $\cF_\real$ & $h(x) = x$ \rule{0mm}{4.5mm} \\
B & $(0,1)$ & $(1,\infty)$  & $w_i \geq 0$; $\sum_{i=1}^k w_i = 1$ & $\cF_\myI^+$ or $\cB^+$ & $h(x) = 1/x$ \rule{0mm}{4.5mm} \\
C & $(0,1)$ & $(-\infty,0)$ & $w_i \geq 0$; $\sum_{i=1}^k w_i > 0$ & $\cF_\myI^+$ or $\cB^+$ & $h(x) = \log x$ \rule{0mm}{4.5mm} \\
D & $(0,1)$ & $\real$       & $w_i \geq 0$; $\sum_{i=1}^k w_i > 0$ & $\cF_\myI^+$ or $\cB^+$ & $h(x) = \Phi^{-1}(x)$ \rule{0mm}{4.5mm} \\
\hline
\end{tabular}

\medskip

\end{table} 

Table \ref{tab:h} shows conditions on the weights, $w_1, \ldots, w_k$,
along with instances of classes $\cF$, so that the generalized linear
combination formula (\ref{eq:GLP}) maps $\cF^k$ into $\cF$.  In the
first type, the link function is defined on the closed unit interval,
and the combination formula operates on the full class $\cF_\real$,
with the traditional linear pool, for which $h(x) = x$ is the identity
function, being the most prominent example.  In the remaining types,
the link function is defined on the open unit interval only, and we
need to restrict attention to $\cF_\myI^+$ or $\cB^+$, with the
harmonic pool and the geometric pool being key examples, occuring when
$h(x) = 1/x$ and $h(x) = \log x$, respectively.  While not being
exhaustive, the listing in the table is comprehensive, in that most
link functions can be adapted to fit one of the types considered.

The following result applies to generalized linear combination formulas 
with nonnegative weights that sum to at most 1.  

\begin{theorem}  \label{th:GLP} 
Suppose that\/ $k \geq 2$, and let\/ $\myI \subseteq \real$ be an
interval.  Consider a generalized linear combination formula\/ $G$ of
the form\/ (\ref{eq:GLP}) with a continuous and strictly monotone link
function, $h$, and weights\/ $w_1, \ldots, w_k$ that are strictly
positive and sum to at most\/ $1$.
\begin{itemize} 
\item[(a)] 
If the link function is defined on the closed unit interval, then\/
$G$ fails to be coherent relative to the class\/ $\cC_\myI$.
\item[(b)] 
If the link function is defined on the open unit interval and
square-integrable, then\/ $G$ fails to be coherent relative to the
class\/ $\cC_\myI^+$.
\end{itemize} 
\end{theorem}

{\em Proof.}  We proceed similarly to the proofs of Theorems
\ref{th:linear} and \ref{th:linear.incoherent}, which correspond to
the special case where $h$ is the identity function and the weights
sum to 1.  As regards part (a) suppose, for a contradiction, that the
generalized linear combination formula $G$ is coherent relative to the
class $\cC_\real$.  Consider the prediction space $(\Omega, \cA,
\myQ)$ with sub-$\sigma$-algebras $\cA_1, \ldots, \cA_k \subseteq
\cA$, CDF-valued random quantities $F_1, \ldots, F_k \in \cC_\real$,
and observation $Y$ such that conditions (C1), (C2) and (C3) of
Definition \ref{def:coherent} hold.  Let $X = h(G(Y))$ and $X_i =
h(F_i(Y))$ for $i = 1, \ldots, k$.  By (C1), (C3), the continuity of
$F_1, \ldots, F_k$, $G$ and $h$, and Theorem \ref{th:ideal}, the
random variables $X$ and $X_1, \ldots X_k$ are identical in
distribution and bounded, and so they share a common finite value of
the variance.  Furthermore, condition (C2) implies that if $i \not= j$
then $X_i \not= X_j$ with positive probability, whence $X_i$ and $X_j$
cannot be linearly dependent.  Consequently,
\[
\var(X) = \sum_{i=1}^k \sum_{j=1}^k w_i w_j \, \cov(X_i,X_j)
        < \sum_{i=1}^k w_i \sum_{j=1}^k w_j \! \left( \max_{1 \leq i \leq k} \var(X_i) \! \right)
        \leq \var(X),
\]
for the desired contradiction.  In part (b) we argue identically,
noting that even though $X$ and $X_1, \ldots, X_k$ may now be
unbounded, they still share a common finite value of the variance, by
the assumption of square-integrability for the link function.  \done

\bigskip
Dawid et al.~(1995) proved in various special cases that generalized
linear combination formulas with nonnegative weights summing to 1 fail
to be coherent relative to the class $\cB^+$ of the nondegenerate
Bernoulli measures.  Furthermore, Dawid et al.~(1985, pp.~282--283)
conjectured that such results hold in broad generality.  While Theorem
\ref{th:GLP} does not apply to Bernoulli measures, as its proof
depends on the continuity of $F_1, \ldots, F_k$ and $G$, we share this
belief.

However, if we allow for individual weights that exceed 1, the
situation may change, as exemplified now in the case of Bernoulli
measures.

\begin{example}  \label{ex:RG2} 
We revisit Example \ref{ex:RG1}, where we identify the ideal
predictive distributions $F_1 = \cL( Y \hsp\hsp | \hsp\hsp \cA_1 )$
and $F_2 = \cL( Y \hsp\hsp | \hsp\hsp \cA_2 )$ with the corresponding
conditional success probabilities $p_1$ and $p_2$, respectively.  It
is then readily seen that
\[
\myQ( \, Y = 0 \mid p_1, p_2 ) 
= \Phi \! \left( (1 + \sigma_2^2)^{1/2} \hsp \Phi^{-1}(p_1) 
  + (1 + \sigma_1^2)^{1/2} \hsp \Phi^{-1}(p_2) \right) \! .   
\]
Hence, the type D generalized linear combination formula
\[
G(y) = \Phi \! \left( w_1 \hsp \Phi^{-1}(p_1) + w_2 \hsp \Phi^{-1}(p_2) \right) 
\] 
with a normal quantile link function, $\Phi^{-1}$, and weights $w_1 >
1$ and $w_2 > 1$ is coherent relative to the class $\cB^+$ of the
nondegenerate Bernoulli measures, or any larger class.
\end{example}

The defining equation (\ref{eq:GLP}) for a generalized linear
combination formula implies that if the weights $w_1, \ldots, w_k$ are
nonnegative and sum to at most 1 then
\[
\var(h(G(Y))) \leq \left( \sum_{i=1}^k w_i \right)^{\! 2} 
\max_{1 \leq i \leq k} \left\{ \var(h(F_i(Y))) \right\} \! .  
\]
In the previous section we applied this type of argument to show that
the traditional linear pool fails to be flexibly dispersive.
Similarly, we conjecture that generalized linear pools with
nonnegative weights that are bounded above fail to be flexibly
dispersive.

It is important to note that the extant literature considers
generalized linear combination formulas that operate on event
probabilities or densities, rather than cumulative distribution
functions.  For example, Dawid et al.~(1995) study geometric and
harmonic linear pools for a binary outcome, while Bj{\o}rnland et
al.~(2011) apply the geometric pool to density forecasts.  Whether or
not the resulting combination formulas and aggregation methods are
coherent and flexibly dispersive, respectively, remains to be
investigated.

\subsection{Spread-adjusted linear pool}  \label{sec:SLP}

The aforementioned limitations of linear and generalized linear pools
suggest that we consider more flexible, nonlinear aggregation methods.
In this section, we focus on the class $\cD_\real^+$, so that we may
identify the cumulative distribution functions $F_1, \ldots, F_k$ with
the corresponding Lebesgue densities $f_1, \ldots, f_k$.

In the context of probabilistic weather forecasts and approximately
neutrally dispersed Gaussian components $f_1, \ldots, f_k$, Berrocal
et al.~(2007), Glahn et al.~(2009) and Kleiber et al.~(2011) observed
empirically that linearly combined predictive distributions are
overdispersed, as confirmed by Theorem \ref{th:linear}.  In an ad hoc
approach, they proposed a nonlinear aggregation method which we now
generalize and refer to as the {\em spread-adjusted linear pool}\/
(SLP), as follows.

To describe this technique, it is convenient to write $F_i(y) =
F_i^0(y - \mu_i)$ and $f_i(y) = f_i^0(y - \mu_i)$, where $\mu_i$ is
the unique median of $F_i \in \cD_\real^+$, for $i = 1, \ldots, k$.
The SLP combined predictive distribution then has cumulative
distribution function and Lebesgue density
\begin{equation}  \label{eq:SLP}
G_c(y) = \sum_{i=1}^k w_i \hsp F_i^0 \! \left( \frac{y-\mu_i}{c} \right) 
\qquad \mbox{and} \qquad 
g_c(y) = \frac{1}{c} \sum_{i=1}^k w_i \hsp f_i^0 \! \left( \frac{y-\mu_i}{c} \right) \! ,
\end{equation}
respectively, where $w_1, \ldots, w_k$ are nonnegative weights that
sum to 1, and $c$ is a strictly positive spread adjustment parameter.
For neutrally dispersed or overdispersed components values of $c < 1$
are appropriate; for example, Table 2 of Berrocal et al.~(2007)
reports estimates ranging from 0.65 to 1.03.  Underdispersed
components may suggest values of $c \geq 1$, and the traditional
linear pool arises when $c = 1$.

We do not know whether or not there are coherent combination formulas
of this type.  However, the SLP method performs well in the
aforementioned applications, and the following result serves to
quantify its flexibility.

\begin{proposition}  \label{th:SLP}
Suppose that\/ $\cL(Y) = F_0 \in \cD_\real^+$ and that\/ $F_1, \ldots,
F_k \in \cD_\real^+$ have medians\/ $\mu_1 \leq \cdots \leq \mu_k$.
Let\/ $Z_c = G_c(Y)$ denote the probability integral transform of the
SLP aggregated predictive cumulative distribution function.  Let\/
$v_0 = 0$ and\/ $p_0 = F_0(\mu_1)$, let\/ $v_i = \sum_{j=1}^i w_j$
and\/ $p_i = F_0(\mu_{i+1}) - F_0(\mu_i)$ for\/ $i = 1, \ldots, k -
1$, and let\/ $v_k = 1$ and\/ $p_k = 1 - F_0(\mu_k)$.  Then as the
spread adjustment parameter\/ $c > 0$ varies, the variance of\/ $Z_c$
attains any positive value less than
\begin{equation}  \label{eq:s0}
\sum_{i=0}^k \, p_i \! \left( v_i - \sum_{j=0}^k p_j \hsp\hsp v_j \right)^{\! 2}.
\end{equation}
\end{proposition}

{\em Proof.} \ As $c \to 0$, the probability integral transform $Z_c$
converges weakly to the discrete probability measure with mass $p_0,
\ldots, p_k$ at $v_0, \ldots, v_k$, which has variance (\ref{eq:s0}).
As $c \to \infty$, it converges weakly to the Dirac measure in
$\frac{1}{2}$.  In view of $\var(Z_c)$ being a continuous function of
the spread-adjustment parameter $c > 0$, this proves the claim. \done

\bigskip
Our next result views the spread-adjusted linear pool as an
aggregation method with parameter space $\Theta = \Delta_{k-1} \times
\real_+$.  While the SLP approach is sufficiently rich in typical
applications, where the individual predictive distributions are
neutrally dispersed or underdispersed, its flexibility is limited.

\begin{theorem}
The spread-adjusted linear pool fails to be flexibly dispersive
relative to the class\/ $\cD_\real^+$.
\end{theorem}

{\em Proof.}  Let $F_0$ be standard normal, and for $i = 1, \ldots, k$
let $F_i$ be normal with mean $m + \frac{i}{m}$ and variance 1.  As $m
\to \infty$, the probability integral transform of the SLP combined
forecast $G_c$ attains values less than $\frac{1}{2}$ with probability
tending to one, irrespectively of the values of the SLP weights $w_1,
\ldots, w_k$ and the spread adjustment parameter $c$.  Thus, if $m$ is
sufficiently large, the variance of the PIT remains below the critical
value of $\frac{1}{12}$ that corresponds to neutral dispersion.  \done

\bigskip
The SLP combination formula (\ref{eq:SLP}) can be generalized to allow
for distinct spread adjustment parameters for the individual
components.  However, such an extension does not admit neutral
dispersion either, and tends not to be beneficial in applications,
unless the component densities have drastically varying degrees of
dispersion.  The assumption of a common spread adjustment parameter
yields a more parsimonious model and stabilizes the estimation.

\subsection{Beta-transformed linear pool}  \label{sec:BLP}

The {\em beta transformed linear pool}\/ (BLP) composites the
traditional linear pool with a beta transform.  Introduced by Ranjan
and Gneiting (2010) in the context of probability forecasts for a
binary event, it generalizes readily to the full class $\cF_\real$ of
the cumulative distribution functions on $\real$.  Specifically, the
BLP combination formula maps $F_1, \ldots, F_k \in \cF_\real$ to
$G_{\alpha, \hsp \beta} \in \cF_\real$, where
\begin{equation}  \label{eq:BLP}
G_{\alpha, \hsp \beta}(y) 
= B_{\alpha, \hsp \beta} \! \left( \sum_{i=1}^k w_i \hsp F_i(y) \! \right) 
\end{equation}
for $y \in \real$.  Here, $w_1, \ldots, w_k$ are nonnegative weights
that sum to $1$, and $B_{\alpha, \hsp \beta}$ denotes the cumulative
distribution function of the beta density with parameters $\alpha > 0$
and $\beta > 0$.  In contrast to the spread-adjusted linear pool, the
value of the BLP aggregated predictive cumulative distribution
function $G_{\alpha, \hsp \beta}$ at $y \in \real$ depends on $F_1,
\ldots, F_k$ only through the values $F_1(y), \ldots, F_k(y)$, in a
locality characteristic that resembles the strong setwise function
property of McConway (1981).  If $F_i$ has Lebesgue density $f_i$ for
$i = 1, \ldots, k$, the aggregated cumulative distribution function
$G_{\alpha, \hsp \beta}$ is absolutely continuous with Lebesgue
density
\[
g_{\alpha, \hsp \beta}(y) 
= \left( \sum_{i=1}^k w_i \hsp\hsp f_i(y) \right)
         b_{\alpha, \hsp \beta} \! \left( \sum_{i=1}^k w_i \hsp\hsp F_i(y) \right) \! ,
\]
where $b_{\alpha, \hsp \beta}$ denotes the beta density with
parameters $\alpha > 0$ and $\beta > 0$.  This nests the traditional
linear pool that arises when $\alpha = \beta = 1$.

The following result concerns the flexibility of the BLP combination
formula (\ref{eq:BLP}) when the cumulative distribution functions $F_0
\in \cC_\real$ and $F_1, \ldots, F_k \in \cC_\real$ are continuous and
the weights $w_1, \ldots, w_k \geq 0$ are fixed, while the
transformation parameters vary.

\begin{proposition}  \label{th:BLP}
Let\/ $Y$ have distribution\/ $F_0 \in \cC_\real$ and suppose that\/
$F_1, \ldots, F_k \in \cC_\real$ are such that
\begin{equation}  \label{eq:support}
{\rm supp}(F_1) \cup \cdots \cup {\rm supp}(F_k) = {\rm supp}(F_0). 
\end{equation}
Let\/ $Z_{\alpha, \hsp \beta} = G_{\alpha, \hsp \beta}(Y)$ denote the
probability integral transform of the BLP aggregated predictice
cumulative distribution function, where the weights are fixed at
strictly positive values that sum to\/ 1.  Then as the transformation
parameters\/ $\alpha > 0$ and\/ $\beta > 0$ vary, the variance of
$Z_{\alpha, \hsp \beta}$ attains any value in the open interval
$(0,\frac{1}{4})$.
\end{proposition}

{\em Proof.} \ The variance of $Z_{\alpha, \hsp \beta}$ depends
continuously on the transformation parameters $\alpha > 0$ and $\beta
> 0$, with $Z_{\alpha, \hsp \alpha}$ converging weakly to the Dirac
measure in $\frac{1}{2}$ as $\alpha \to \infty$, so that
$\var(Z_{\alpha, \hsp \alpha}) \to 0$ as $\alpha \to \infty$.  If we
can demonstrate the existence of a sequence $(\alpha, \beta(\alpha))
\to (0,0)$ such that $G_{\alpha, \hsp \beta(\alpha)}(y_0) =
\frac{1}{2}$, where $y_0$ is any median of $F_0$, the proof is
complete, as the corresponding probability integral transform
$Z_{\alpha, \hsp \beta(\alpha)}$ converges weakly to the Bernoulli
measure with success probability $\frac{1}{2}$, so that
$\var(Z_{\alpha, \hsp \beta(\alpha)}) \to \frac{1}{4}$ as $\alpha \to
0$.

We thus strive to find a sequence $(\alpha, \beta(\alpha)) \to (0,0)$
such that $G_{\alpha, \hsp \beta(\alpha)}(y_0) = B_{\alpha, \hsp
  \beta(\alpha)}(u_0) = \frac{1}{2}$, where $u_0 = \sum_{i=1}^{k} w_i
\hsp F_i(y_0) \in (0,1)$ by the support condition (\ref{eq:support}).
First we show that for every $\alpha > 0$ there exists a unique
$\beta(\alpha) > 0$ such that $B_{\alpha, \, \beta(\alpha)}(u_0) =
\frac{1}{2}$; then we prove that $\beta(\alpha) \to 0$ as $\alpha \to
0$.  As regards the first claim, three cases are to be distinguished.
If $u_0 < \frac{1}{2}$ then $B_{\alpha, \hsp \alpha}(u_0) <
\frac{1}{2}$ and $B_{\alpha, \hsp \beta}(u_0) \to 1$ as $\beta \to
\infty$, and continuity and monotonicity with respect to $\beta$ imply
the existence of a unique $\beta(\alpha) > \alpha$ such that
$B_{\alpha, \hsp \beta(\alpha)}(u_0) = \frac{1}{2}$.  If $u_0 =
\frac{1}{2}$ the choice $\beta(\alpha) = \alpha$ is unique.  If $u_0 >
\frac{1}{2}$ then $B_{\alpha, \hsp \alpha}(u_0) > \frac{1}{2}$ and
$B_{\alpha, \hsp \beta}(u_0) \to 0$ monotonically as $\beta \to 0$,
and thus there exists a unique $\beta(\alpha) < \alpha$ such that
$B_{\alpha, \hsp \beta(\alpha)}(u_0) = \frac{1}{2}$.  To prove the
second claim, suppose that $\beta(\alpha) > \beta_0 > 0$ for a
sequence $\alpha \to 0$.  Then as $\alpha \to 0$ the beta distribution
with parameters $\alpha$ and $\beta(\alpha)$ has mean $\alpha/(\alpha
+ \beta(\alpha)) \to 0$, whereas its median $u_0$ remains fixed and
strictly positive, for the desired contradiction.  \done

\bigskip
The next result views the beta-transformed linear pool as an
aggregation method with parameter space $\Theta = \Delta_{k-1} \times
\real_+^2$.  

\begin{theorem}
The beta transformed linear pool is exchangeably flexibly dispersive
relative to the class\/ $\cC_\myI^+$, for every interval\/ $I
\subseteq \real$.
\end{theorem}

{\em Proof.} \ If $F_0 \in \cC_\myI^+$ and $F_1, \ldots, F_k \in
\cC_\myI^+$, the support condition (\ref{eq:support}) is satisfied.
We may thus apply Theorem \ref{th:BLP} with weights $w_1 = \cdots =
w_k = \frac{1}{k}$.  \done

\bigskip
As hinted at in Section \ref{sec:aggregation}, the spread-adjusted and
beta transformed linear pools can be applied in the case $k = 1$ of a
single probabilistic forecast to provide calibration and dispersion
adjustments.  If the original predictive distribution is symmetric
about its center, the symmetry is retained by the spead-adjusted
linear pool, and broken by the beta transformed linear pool, except
when $\alpha = \beta$.

In practice, the BLP weights $w_1, \ldots, w_k$ and transformation
parameters $\alpha, \beta > 0$ are estimated from training data, say
$\{ (F_{1j}, \ldots, F_{kj}, y_j) : j = 1, \ldots, J \}$.  If the
predictive cumulative distribution functions $F_{1j}, \ldots, F_{kj}$
are absolutely continuous with Lebesgue densities $f_{1j}, \ldots,
f_{kj}$ for $j = 1, \ldots, J$, the aggregated predictive
distributions are also absolutely continuous, and our preferred
estimation technique is to maximize the mean (or sum) of the
logarithmic score (Matheson and Winkler 1976; Gneiting and Raftery
2007) over the training data, namely
\begin{eqnarray}  \label{eq:log.likelihood}
\ell \hsp (w_1, \ldots, w_k; \alpha, \beta)
& = & \sum_{j=1}^J \log \hsp (g_{\alpha,\beta}(y_j)) \\ \nonumber
& = & \sum_{j=1}^J \log \! \left( \sum_{i=1}^k w_i \hsp f_{ij}(y_j) \right)
      \; + \; \sum_{j=1}^J \log \! \left( b_{\alpha,\beta} \! \left(\sum_{i=1}^k
      w_i \hsp F_{ij}(y_j) \right) \right) \\ \nonumber
& = & \sum_{j=1}^J \left( (\alpha-1) \log \! \left( \sum_{i=1}^k
      w_i \hsp F_{ij}(y_j) \right) + (\beta-1) \log \! \left( \! 1-\sum_{i=1}^k
      w_i \hsp F_{ij}(y_j) \right) \right) \\ \nonumber
&   & + \;\, \sum_{j=1}^J \log \! \left( \sum_{i=1}^k w_i \hsp f_{ij}(y_j) \right)
      \, - \; J \log \hsp {\rm B}(\alpha,\beta),
\end{eqnarray}
where ${\rm B}$ denotes the classical beta function.  The logarithmic
score is simply the logarithm of the value that the density forecast
attains at the realizing observation.  It is positively oriented, that
is, the higher the score, the better, and it is proper, in the sense
that truth telling is an expectation maximizing strategy.
Alternatively, the corresponding estimates can be viewed as maximum
likelihood estimates under the assumption of independence between the
training cases.

The optimization can be carried out numerically using the method of
scoring, for which we give details in Appendix B.  Approximate
standard errors for the estimates can be obtained in the usual way, by
evaluating and inverting the Hessian matrix for the mean logarithmic
score or log likelihood function.  However, the estimates of the
weights $w_1, \ldots, w_k$ need to be nonnegative.  Thus, if
unconstrained optimization results in negative weights, we turn to the
active barrier algorithm implemented in the constrained optimization
routine {\sc constrOptim} in {\sc R} ({\sc{R}} Development Core Team
2011).  Similarly, linear, generalized linear and spread-adjusted
linear combination formulas can be fitted by maximizing the mean
logarithmic score over training data.  For reasons of simplicity and
tradition, we frequently refer to the resulting estimates as maximum
likelihood estimates.

\section{Simulation and data examples}  \label{sec:examples}

We now illustrate and complement our theoretical results in simulation
and data examples on density forecasts.  This corresponds to the
prediction space setting, where the CDF-valued random quantities $F_1,
\ldots, F_k$ are absolutely continuous almost surely, and thus can be
identified with random Lebesgue densities $f_1, \ldots, f_k$.
Throughout the section, we fit combination formulas by maximizing the
mean logarithmic score over training data, in the ways described above
and in Appendix B.  To lighten the notation, we use the acronyms PIT,
TLP, SLP and BLP to refer to the probability integral transform and
the traditional, spread-adjusted and beta transformed linear pool,
respectively.

The recent work of Ranjan and Gneiting (2010) and Clements and Harvey
(2011) contains a wealth of simulation and data examples on the
combination of probability forecasts for a binary event.  In the
concluding Section \ref{sec:discussion} we summarize these experiences
and relate them to the findings in the case studies hereinafter.

\subsection{Simulation example}  \label{sec:sim}

\begin{table}[p]
\centering
\caption{Maximum likelihood estimates with approximate standard errors
  (in brackets) for the parameters of the combined density forecasts
  in the simulation example. \label{tab:sim.MLE}} \footnotesize
\medskip
\begin{tabular}{lcccccc}
\hline
\hline
    & $w_1$ & $w_2$ & $w_3$ & $c$ & $\alpha$ & $\beta$ \\
\hline
TLP & 0.212 (0.083) & 0.254 (0.084) & 0.534 (0.080) & ---           & ---           & --- \\
SLP & 0.257 (0.060) & 0.283 (0.061) & 0.460 (0.059) & 0.783 (0.030) & ---           & --- \\
BLP & 0.256 (0.057) & 0.293 (0.057) & 0.451 (0.054) & ---           & 1.492 (0.062) & 1.440 (0.059) \\
\hline
\end{tabular}
\smallskip
\end{table}

\begin{table}[p]
\centering
\caption{Variance of the PIT (dispersion) and root mean variance of
  the density forecast (sharpness) in the simulation example, for the
  test set. A value of $\frac{1}{12}$ or about $0.083$ for the
  variance of the PIT indicates neutral
  dispersion.  \label{tab:sim.dispersion.sharpness}} \footnotesize
\medskip
\begin{tabular}{lcc}
\hline
\hline
      & var(PIT) & RMV \\
\hline
$f_1$ & $0.081$ & $1.79$ \\
$f_2$ & $0.086$ & $1.79$ \\
$f_3$ & $0.085$ & $1.73$ \\
TLP   & $0.066$ & $1.94$ \\
SLP   & $0.081$ & $1.62$ \\
BLP   & $0.084$ & $1.57$ \\
\hline
\end{tabular}
\smallskip
\end{table}

\begin{figure}[p]
\centering
\includegraphics[width = 0.75\textwidth]{\mypath 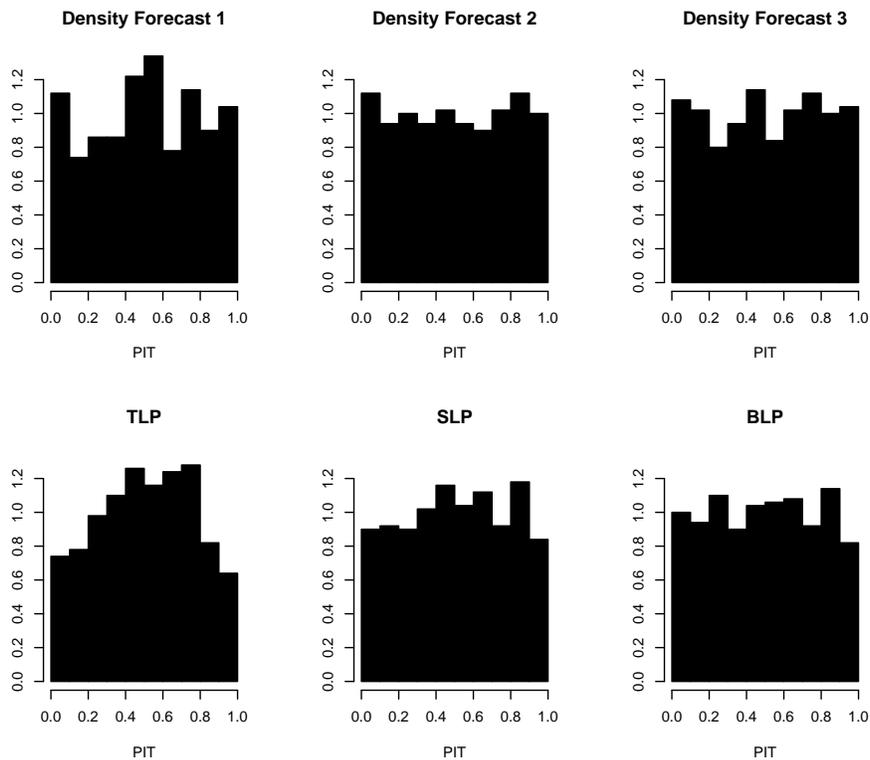}
\caption{PIT histograms for the individual and combined density
  forecasts in the simulation example, for the test
  set. \label{fig:sim.PIT}}
\end{figure}

In this simulation example, the data generating process for the
observation, $Y$, is the regression model
 \begin{equation}  \label{eq:DGP}
Y = X_0 + a_1 X_1 + a_2 X_2 + a_3 \hsp X_3 + \epsilon,
\end{equation}
where $a_1, a_2$ and $a_3$ are real constants, and $X_0, X_1, X_2,
X_3$ and $\epsilon$ are independent, standard normal random variables.
The individual predictive distributions rest on partial knowledge of
the data generating process, in that density forecast $f_1$ has access
to the covariates $X_0$ and $X_1$, but not to $X_2$ or $X_3$, and
similarly for $f_2$ and $f_3$.  Thus, we seek to combine the density
forecasts
\[
f_1 = \cN(X_0 + a_1 X_1, 1 + a_2^2 + a_3^2),
\]
\[
f_2 = \cN(X_0 + a_2 X_2, 1 + a_1^2 + a_3^2) \quad \mbox{and} \quad
f_3 = \cN(X_0 + a_3 \hsp X_3, 1 + a_1^2 + a_2^2),
\]
where $X_0$ stands for shared, public information, while $X_1$, $X_2$
and $X_3$ represent proprietary information sets.  The density
forecasts represent the true conditional distributions under the
regression model (\ref{eq:DGP}), given the corresponding partial
information, as represented by the $\sigma$-algebras $\cA_1 =
\sigma(X_0, X_1)$, $\cA_2 = \sigma(X_0, X_2)$ and $\cA_3 = \sigma(X_0,
X_3)$, respectively.  Hence, the forecasts are ideal in the sense of
Definition \ref{def:ideal}, and by Theorem \ref{th:ideal} they are
both probabilistically calibrated and marginally calibrated.  

We estimate the TLP, SLP and BLP combination formulas on a simple
random sample $\{ (f_{1j}, f_{2j}, f_{3j}, Y_j) : j = 1, \ldots, J \}$
of size $J = 500$ from the joint distribution of the forecasts and the
observation, and evaluate on an independent test sample of the same
size.  The regression coefficients in the data generating model
(\ref{eq:DGP}) are taken to be $a_1 = a_2 = 1$ and $a_3 = 1.1$, so
that $f_3$ is a more concentrated, sharper density forecast than $f_1$
and $f_2$.

\begin{table}[t]
\centering
\caption{Mean logarithmic score for the individual and combined
  density forecasts in the simulation example, for the training set
  and the test set. \label{tab:sim.score}} \footnotesize
\bigskip
\begin{tabular}{lcc}
\hline
\hline
      & Training & Test \\
\hline
$f_1$ & $-2.025$ & $-2.018$ \\
$f_2$ & $-2.017$ & $-2.022$ \\
$f_3$ & $-1.956$ & $-1.992$ \\
TLP   & $-1.907$ & $-1.922$ \\
SLP   & $-1.871$ & $-1.892$ \\
BLP   & $-1.865$ & $-1.886$ \\
\hline
\end{tabular}
\medskip
\end{table}

Table \ref{tab:sim.MLE} shows maximum likelihood estimates, along with
approximate standard errors, for TLP, SLP and BLP combination
formulas.  For all three methods, the weight estimate is highest for
$f_3$, whereas the estimates for $f_1$ and $f_2$ are smaller and not
significantly different from each other.  The SLP spread adjustment
parameter $c$ is estimated at 0.78, and the BLP transformation
parameters $\alpha$ and $\beta$ at 1.49 and 1.44, respectively.

The PIT histograms for the various types of density forecasts over the
test set are displayed in Figure \ref{fig:sim.PIT}, with complementary
results shown in Table \ref{tab:sim.dispersion.sharpness}.  In
addition to the variance of the PIT, which is our standard measure of
dispersion, the table quantifies sharpness in terms of the root mean
variance (RMV), that is, the square root of the average of the
variance of the predictive density over the evaluation set.  The
component forecasts $f_1$, $f_2$ and $f_3$ are probabilistically
calibrated and thus show uniform empirical PIT histograms, up to
sample fluctuations.  As mandated by Theorem \ref{th:linear}, the
linearly combined TLP density forecast is overdispersed and lacks
sharpness.  The SLP and BLP agrregated density forecasts show nearly
uniform PIT histograms; they are approximately neutrally dispersed and
much sharper than their competitors.

Table \ref{tab:sim.score} shows the mean logarithmic score for the
various types of density forecasts.  The best individual density
forecast is $f_3$, because it is sharper than $f_1$ and $f_2$.  The
linearly combined density forecast outperforms the individual density
forecasts, even though it is overdispersed.  The nonlinearly
aggregated SLP and BLP density forecasts show higher scores than any
of the individual or linearly combined forecasts, both for the
training data, where this is trivially true, as the nonlinear methods
nest the traditional linear pool, and for the test data, where such
cannot be guaranteed.

\subsection{Density forecasts for daily maximum temperature at
            Seattle-Tacoma Airport} \label{sec:temperature}

With estimates of some one-third of the economy, as well as much of
human activity in general, being weather sensitive (Dutton 2002),
there is a critical need for calibrated and sharp probabilistic
weather forecasts, to allow for optimal decision making under inherent
environmental uncertainty.

In practice, probabilistic weather forecasts rely on ensemble
prediction systems.  An ensemble system comprises multiple runs of a
numerical weather prediction model, with the runs differing in the
initial conditions and/or the details of the mathematical
representation of the atmosphere (Palmer 2002; Gneiting and Raftery
2005).  Here we consider two-days ahead forecasts of daily maximum
temperature at Seattle-Tacoma Airport, based on the University of
Washington Mesoscale Ensemble (Eckel and Mass 2005), which employs a
regional numerical weather prediction model over the Pacific
Northwest, with initial and lateral boundary conditions supplied by
eight distinct weather centers.  A brief description of the ensemble
members is given in Table \ref{tab:UWME}.

\begin{table}[p]

\centering

\caption{Composition of the eight-member University of Washington
  Mesoscale Ensemble (Eckel and Mass 2005), with member acronyms and
  organizational sources for initial and lateral boundary conditions.
  The United States National Centers for Environmental Prediction
  supply two distinct sets of initial and lateral boundary conditions,
  namely, from its Global Forecast System (GFS) and Limited-Area
  Mesoscale Model (ETA).  \label{tab:UWME}}
\medskip
\footnotesize
\begin{tabular}{llll}
\hline
\hline
Index & Acronym & Source of Initial and Lateral Boundary Conditions \\
\hline
1 & GFS  & National Centers for Environmental Prediction \\
2 & CMCG & Canadian Meteorological Centre \\
3 & ETA  & National Centers for Environmental Prediction \\
4 & GASP & Australian Bureau of Meteorology \\
5 & JMA  & Japanese Meteorological Agency \\
6 & NGPS & Fleet Numerical Meteorology and Oceanography Center \\
7 & TCWB & Taiwan Central Weather Bureau \\
8 & UKMO & United Kingdom Met Office \\
\hline
\end{tabular}
\bigskip

\caption{Maximum likelihood estimates for the predictive standard
  deviation, $\sigma_i$, for the individual, member specific density
  forecasts in the temperature example. \label{tab:temp.sd}}
\medskip
\footnotesize
\begin{tabular}{cccccccc}
\hline \hline
$\sigma_1$ & $\sigma_2$ & $\sigma_3$ & $\sigma_4$ & $\sigma_5$ & $\sigma_6$ & $\sigma_7$ & $\sigma_8$ \\
\hline
1.966 & 2.051 & 2.119 & 2.214 & 1.958 & 2.055 & 2.084 & 1.995 \\
\hline
\end{tabular}
\bigskip

\caption{Maximum likelihood estimates for the parameters of the
  combined density forecasts in the temperature example, including the
  Bayesian model averaging (BMA) approach of Raftery et
  al.~(2005). \label{tab:temp.MLE}} \footnotesize
\medskip
\begin{tabular}{lcccccccccccc}
\hline
\hline
    & $w_1$ & $w_2$ & $w_3$ & $w_4$ & $w_5$ & $w_6$ & $w_7$ & $w_8$ & $c$ & $\alpha$ & $\beta$ & $\sigma$ \\
\hline
TLP & 0.394 & 0.005 & 0.000 & 0.000 & 0.317 & 0.030 & 0.144 & 0.109 & ---   & ---   & ---   & ---   \\
SLP & 0.304 & 0.080 & 0.000 & 0.085 & 0.216 & 0.051 & 0.172 & 0.090 & 0.768 & ---   & ---   & ---   \\
BLP & 0.295 & 0.079 & 0.000 & 0.083 & 0.230 & 0.062 & 0.173 & 0.076 & ---   & 1.467 & 1.467 & ---   \\
BMA & 0.305 & 0.075 & 0.000 & 0.081 & 0.216 & 0.056 & 0.170 & 0.098 & ---   & ---   & ---   & 1.566 \\
 \hline
\end{tabular}
\bigskip

\caption{Variance of the PIT (dispersion) and root mean variance of
  the density forecast (sharpness) in the temperature example, for the
  test period.  A value of $\frac{1}{12}$ or about $0.083$ for the
  variance of the PIT indicates neutral
  dispersion. \label{tab:temp.dispersion.sharpness}} \footnotesize
\medskip
\begin{tabular}{lcc}
\hline
\hline
      & var(PIT) & RMV \\
\hline
$f_1$ & $ 0.070$ & $1.97$ \\
$f_2$ & $ 0.067$ & $2.05$ \\
$f_3$ & $ 0.069$ & $2.12$ \\
$f_4$ & $ 0.068$ & $2.21$ \\
$f_5$ & $ 0.070$ & $1.96$ \\
$f_6$ & $ 0.073$ & $2.06$ \\
$f_7$ & $ 0.074$ & $2.08$ \\
$f_8$ & $ 0.069$ & $2.00$ \\
TLP   & $ 0.057$ & $2.15$ \\
SLP   & $ 0.070$ & $1.79$ \\
BLP   & $ 0.072$ & $1.77$ \\
BMA   & $ 0.070$ & $1.80$ \\
\hline
\end{tabular}
\end{table}

Our training period ranges from January 1, 2006 to August 12, 2007,
with a few days missing in the data record, for a total of 500
training cases.  The test period extends from August 13, 2007 to June
30, 2009, for a total of 559 cases.

Each ensemble member is a point forecast, which can be viewed as the
most extreme form of an underdispersed density forecast.  To address
the underdispersion and obtain approximately neutrally dispersed
components, we use the maximum likelihood method on the training data
to fit, for each ensemble member $i = 1, \ldots, 8$ individually, a
Gaussian predictive density of the form
\[
f_i = \cN(a_i + b_i \hsp x_{ij}, \sigma_i^2). 
\]
Here $x_{ij}$ is the point forecast from the $i$\/\/th ensemble member
on day $j$, $a_i$ and $b_i$ are member specific linear bias correction
parameters, and $\sigma_i$ is the member specific predictive standard
deviation.  From Table \ref{tab:temp.sd} we see that the estimates for
$\sigma_1, \ldots, \sigma_8$ range from 1.958 to 2.214.

Next we combine the eight individual density forecasts.  Table
\ref{tab:temp.MLE} shows maximum likelihood estimates for TLP, SLP and
BLP combination formulas.  For all three methods, the GFS member,
$f_1$, obtains the highest weight and the ETA member, $f_3$, the
lowest weight.  This can readily be explained, in that both members
have a common institutional origin, and thus are highly correlated,
whence the more competitive GFS member subsumes the weight of the ETA
member.  The SLP spread adjustment parameter is estimated at 0.768,
and the BLP transformation parameters both at 1.467.

\begin{figure}[t]
\centering
\includegraphics[width = \textwidth]{\mypath 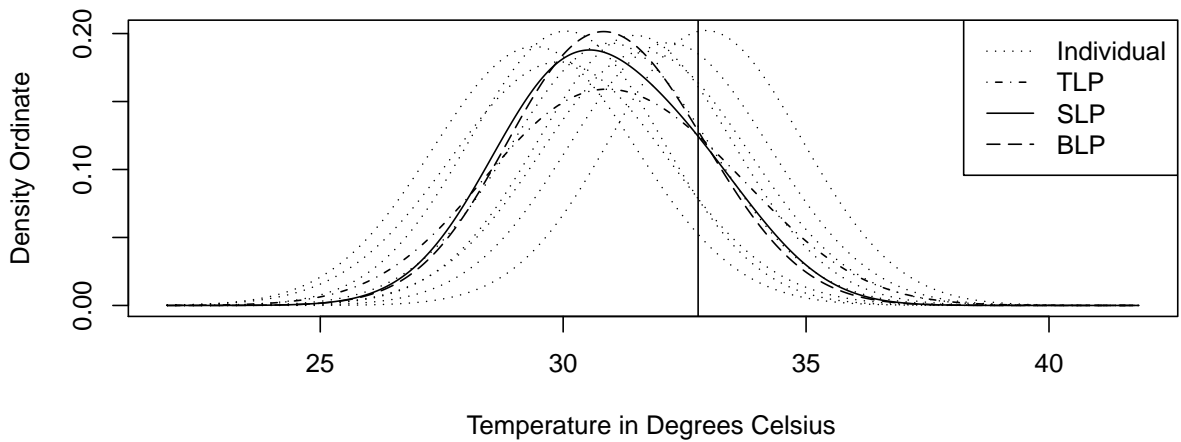}
\caption{Two-day ahead density forecasts for the maximum temperature
  at Seattle-Tacoma Airport on June 28, 2008.  The vertical line is at
  the verifying maximum, at 32.8 degrees Celsius or 91 degrees
  Fahrenheit.\label{fig:temp.density}}
\end{figure}

Figure \ref{fig:temp.density} illustrates the various density
forecasts for June 28, 2008, an unusually hot day at Seattle-Tacoma
Airport with a verifying maximum temperature of 32.8 degrees Celsius
or 91 degrees Fahrenheit.  The member specific individual density
forecasts are shown by the dotted lines, and the linearly combined TLP
forecast by the dash-dotted line.  The nonlinearly aggregated SLP and
BLP density forecasts, which are shown by the solid and dashed line,
respectively, are sharper than the TLP density.

\begin{figure}[p]
\centering
\bigskip
\includegraphics[width = \textwidth]{\mypath 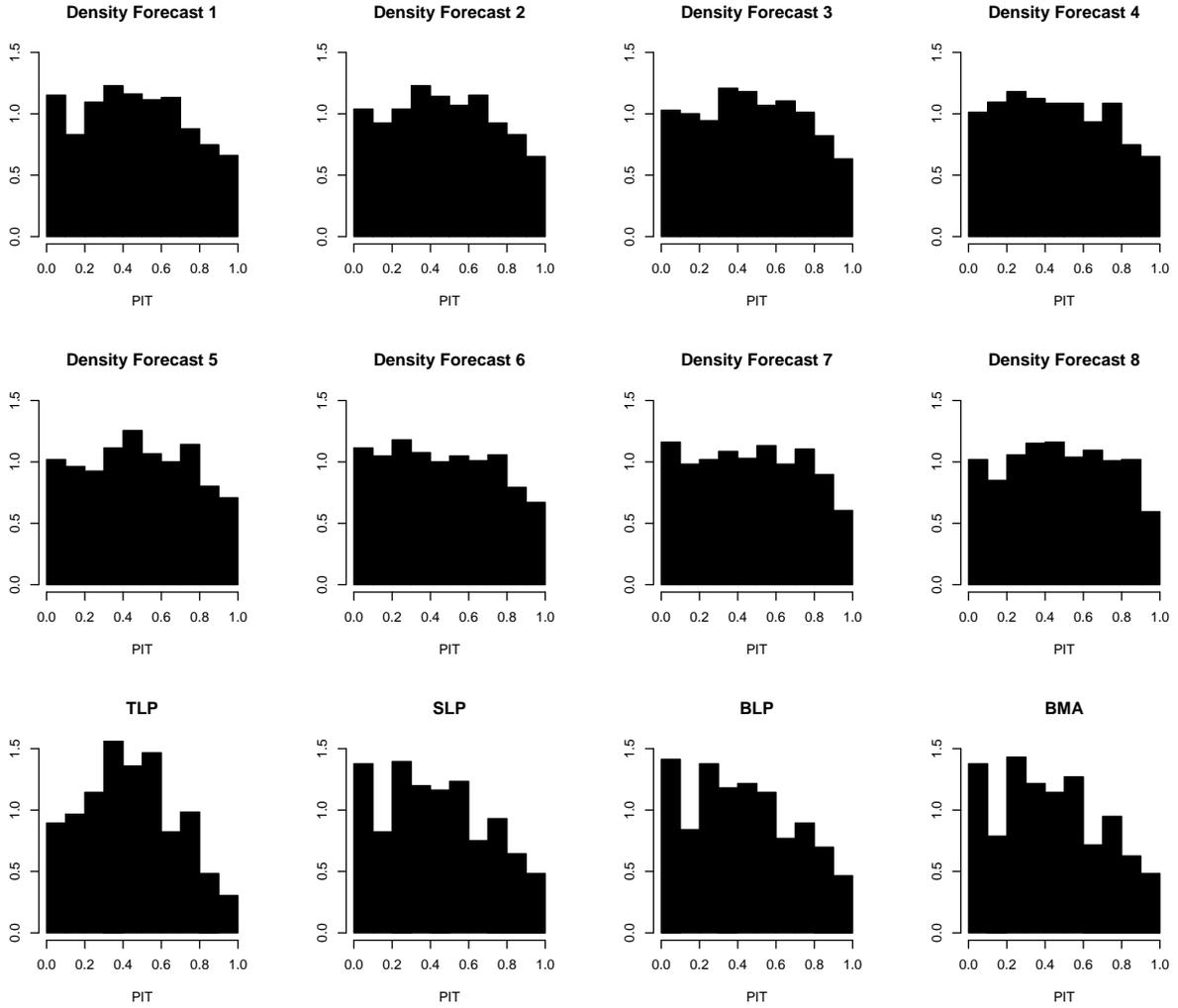}
\caption{PIT histograms for the individual and combined density forecasts in
  the temperature example, for the test period. \label{fig:temp.PIT}}
\end{figure}

PIT histograms for the test period are shown in Figure
\ref{fig:temp.PIT}, along with summary measures of dispersion and
sharpness in Table \ref{tab:temp.dispersion.sharpness}.  The
individual, member specific density forecasts tend to be a bit
overdispersed.  The linearly aggregated TLP density forecast is much
more severely overdispersed, as reflected by an inverse U-shaped and
skewed PIT histogram.  Of course, the overdispersion is not
surprising, as it is a direct consequence of Theorem \ref{th:linear}.
The SLP and BLP aggregated density forecasts show somewhat rough and
skewed, yet more nearly uniform PIT histograms.

\begin{table}[t]
\centering
\caption{Mean logarithmic score for the individual and combined
  density forecasts in the temperature example, for the training
  period and the test period. \label{tab:temp.score}}
\footnotesize
\medskip
\begin{tabular}{lccc}
\hline
\hline
      & Training & Test \\
\hline
$f_1$ & $-2.091$ & $-2.088$ \\
$f_2$ & $-2.134$ & $-2.071$ \\
$f_3$ & $-2.167$ & $-2.093$ \\
$f_4$ & $-2.211$ & $-2.172$ \\
$f_5$ & $-2.088$ & $-2.043$ \\
$f_6$ & $-2.136$ & $-2.143$ \\
$f_7$ & $-2.150$ & $-2.131$ \\
$f_8$ & $-2.107$ & $-2.041$ \\
TLP   & $-2.027$ & $-2.010$ \\
SLP   & $-1.990$ & $-1.961$ \\
BLP   & $-1.988$ & $-1.960$ \\
BMA   & $-1.992$ & $-1.963$ \\
\hline
\end{tabular}
\bigskip
\end{table}

These results are corroborated by Table \ref{tab:temp.score}, which
shows the mean logarithmic score for the various types of density
forecasts, both for the training period and the test period.  The
linearly combined TLP forecast shows a higher score than any of the
individual density forecasts, which attests to the benefits of
aggregation.  Nevertheless, the linearly combined density forecast is
suboptimal, because it is overdispersed and lacks sharpness, and thus
it is outperformed by the nonlinearly aggregated SLP and BLP density
forecasts.

Finally, we compare to the Bayesian model averaging (BMA; Raftery et
al.~2005) technique, which is a state of the art approach to
generating density forecasts from forecast ensembles.  The BMA density
forecast for day $j$ is of the form
\begin{equation}  \label{eq:BMA}
g = \sum_{i=1}^8 w_i \; {\cal N} (a_i + b_i \hsp x_{ij}, \, \sigma^2),
\end{equation}
with BMA weights, $w_1, \ldots, w_8$, that are nonnegative and sum to
1, member specific bias parameters $a_i$ and $b_i$ for $i = 1, \ldots,
8$, and a common variance parameter, $\sigma^2$.  In view of our
individual density forecasts being Gaussian, the TLP and BMA densities
are of the same functional form.  However, there is a conceptual
difference, in that the TLP weights are fitted conditionally on the
individual density forecasts.  Thus, a two-stage procedure is used, in
which the member specific component densities are estimated first, and
only then the weights, with the components held fixed.  In contrast,
the BMA method estimates the weights, $w_1, \ldots, w_8$, and the
common spread parameter, $\sigma$, for the component forecasts in the
Gaussian mixture model (\ref{eq:BMA}) simultaneously.  While the BMA
method can be employed with member specific spread parameters, the
assumption of a common spread parameter stabilizes the estimation
algorithm and does not appreciably deteriorate the predictive
performance (Raftery et al.~2005).

Table \ref{tab:temp.MLE} shows maximum likelihood estimates for the
BMA parameters, obtained with the {\sc R} package {\sc ensembleBMA}
(Fraley et al.~2011).  The BMA weights echo the SLP weights.  The BMA
spread parameter $\sigma$ is estimated at 1.566 and differs from the
predictive standard deviations for the member specific density
forecasts in Table \ref{tab:temp.sd} by factors ranging from 0.707 to
0.800, much in line with our estimate of 0.768 for the SLP spread
adjustment parameter, $c$.  Thus, the SLP and BMA density forecasts
are very much alike, which is confirmed by the PIT histograms in
Figure \ref{fig:temp.PIT}, the summary measures in Table
\ref{tab:temp.dispersion.sharpness} and the logarithmic scores in
Table \ref{tab:temp.score}.  In Figure \ref{fig:temp.density} the
graphs for the SLP and BMA density forecasts are nearly identical and
lie essentially on top of each other, and so we refrain from plotting
the BMA density.

\subsection{Density forecasts for S\&P 500 returns}  \label{sec:SP}

In this final data example, we follow Diebold et al.~(1998) in
considering S\&P 500 log returns for the period of July 3, 1962 to
December 29, 1995.  The data record through December 1978 is used as
training set, for a total of 4,133 training cases.  All estimates
reported are maximum likelihood fits on the training period obtained
with the {\sc R} package {\sc fGarch} (Wuertz and Rmetrics Core Team
2007).  The balance of the record is used as test period, for a total
of 4,298 one-day ahead density forecasts.

The first component forecast, $f_1$, is based on a generalized
autoregressive conditional heteroscedasticity (GARCH; Bollerslev 1986)
specification for the variance structure.  With $r_t$ denoting the log
return on day $t$, our GARCH(1,1) model assumes that $r_t = \sigma_t
\hsp \epsilon_{\hsp t}$, where $\epsilon_{\hsp t}$ is Student-$t$
distributed with $\nu$ degrees of freedom and variance 1, while
$\sigma_t$ evolves dynamically as
\[
\sigma_t^2 =  \omega + \alpha \hsp r_{t-1}^2 + \beta \hsp \sigma_{t-1}^2.
\]
The maximum likelihood estimates for the GARCH parameters are $\omega
= 0.000$, $\alpha = 0.089$, $\beta = 0.903$ and $\nu = 9.25$.

The second component forecast, $f_2$, is based on a standard moving
average (MA) model for the mean dynamics, which assumes that $r_t =
Z_t + \theta Z_{t-1}$, where $\{ Z_t \}$ is a Gaussian white noise
process with mean zero and variance $\sigma^2$.  The maximum
likelihood estimates for the MA parameters are $\theta = 0.252$ and
$\sigma = 0.00736$.

\begin{figure}[t]
\centering
\includegraphics[width = \textwidth]{\mypath 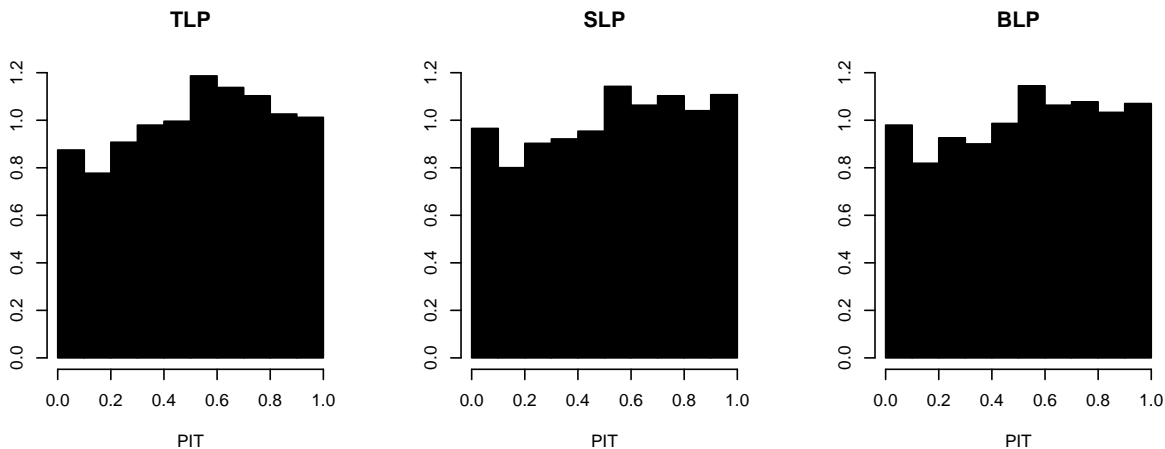}
\caption{PIT histograms for the combined density forecasts in the S\&P
  500 example, for the test period. \label{fig:SP500.PIT}}
\bigskip
\end{figure}

Our goal now is to combine the density forecasts $f_1$ and $f_2$.
Table \ref{tab:SP500.MLE} shows maximum likelihood estimates for TLP,
SLP and BLP combination formulas.  For all three methods, the
conditionally heteroscedastic density forecast $f_1$ obtains a much
higher weight than the simplistic density forecast $f_2$.  The SLP
spread adjustment parameter is estimated at 0.940, and the BLP
transformation parameters $\alpha$ and $\beta$ at 1.100 and 1.081.
This suggests that the overdispersion of the TLP density forecast is
quite mild, which is confirmed by the PIT histogram in Figure
\ref{fig:SP500.PIT} and the summary measures in Table
\ref{tab:SP500.dispersion.sharpness}.

Table \ref{tab:SP500.score} shows the mean logarithmic score for the
various types of probabilistic forecasts.  The TLP density forecast
performs slightly better than the individual component $f_1$, with a
score that is very slightly lower than for the nonlinearly aggregated
SLP and BLP density forecasts, both for the training and the test
period.  As observed by Geweke and Amisano (2011), there is little
reward for using more elaborate, less parsimonious aggregation methods
for density forecasts of S\&P 500 returns.\footnote{The logarithmic
  scores reported by Geweke and Amisano (2011) are summed, rather than
  averaged, and apply to percent log returns, rather than log returns.
  Adjusted for these differences, they are comparable to the scores
  in Table \ref{tab:SP500.score}.}

\begin{table}[t]
\centering
\caption{Maximum likelihood estimates of the parameters for the
  combined density forecasts in the S\&P 500
  example. \label{tab:SP500.MLE}} \footnotesize
\bigskip
\begin{tabular}{lccccc}
\hline
\hline
    & $w_1$ & $w_2$ & $c$ & $\alpha$ & $\beta$ \\
\hline
TLP & 0.821 & 0.179 & ---   & ---   & ---   \\
SLP & 0.756 & 0.244 & 0.940 & ---   & ---   \\
BLP & 0.758 & 0.242 & ---   & 1.100 & 1.081 \\
\hline
\end{tabular}
\medskip
\end{table}

\begin{table}[t]
\centering
\caption{Variance of the PIT (dispersion) and root mean variance of
  the density forecast (sharpness) in the S\&P 500 example, for the
  test period.  A value of $\frac{1}{12}$ or about $0.083$ for the
  variance of the PIT indicates neutral
  dispersion. \label{tab:SP500.dispersion.sharpness}} \footnotesize
\bigskip
\begin{tabular}{lcc}
\hline
\hline
      & var(PIT) & RMV \\
\hline
$f_1$ & $ 0.081$ & $9.23 \times 10^{-3}$ \\
$f_2$ & $ 0.078$ & $7.36 \times 10^{-3}$ \\
TLP   & $ 0.079$ & $8.89 \times 10^{-3}$ \\
SLP   & $ 0.083$ & $8.22 \times 10^{-3}$ \\
BLP   & $ 0.084$ & $8.22 \times 10^{-3}$ \\
\hline
\end{tabular}
\medskip
\end{table}

\begin{table}[t]
\centering
\caption{Mean logarithmic score for the individual and combined
  density forecasts in the S\&P 500 example, for the training period
  and the test period. \label{tab:SP500.score}}
\bigskip
\footnotesize
\begin{tabular}{lcc}
\hline
\hline
      & Training & Test \\
\hline
$f_1$ & 3.606 & 3.458 \\
$f_2$ & 3.492 & 3.247 \\
TLP   & 3.612 & 3.469 \\
SLP   & 3.614 & 3.470 \\
BLP   & 3.614 & 3.470 \\
\hline
\end{tabular}
\bigskip
\end{table}

Finally, we consider the predictive performance of a more 
comprehensive predictive model, which addresses both the first and the
second order dynamics, in that $r_t = \mu_t + \epsilon_t$ where $\{
\mu_t \}$ and $\{ \epsilon_t \}$ are MA(1) and Student-$t$ GARCH(1,1)
processes, respectively.  The maximum likelihood estimates in this
mixed specification are $\theta = 0.269$ and $\sigma = 0.00736$ for
the MA parameters, and $\omega = 0.000$, $\alpha = 0.098$, $\beta =
0.892$ and $\nu = 8.284$ for the GARCH parameters.  The resulting
density forecast can be thought of as combining information sets with
respect to the first and second order dynamics, as opposed to
combining the corresponding component forecasts $f_1$ and $f_2$.  It
outperforms the other types of density forecasts and achieves a mean
logarithmic score of 3.638 for the training period and 3.473 for the
test period.

\section{Discussion} \label{sec:discussion}

We have studied methods for combining predictive distributions.  From
a theoretical perspective, our approach deviates from previous work in
major ways.  Technically, we operate in terms of prediction spaces and
cumulative distribution functions, which allows for a unified
treatment of all real-valued predictands including, for example,
density forecasts for continuous variables and probability forecasts
of a binary event.  Conceptually, our work is motivated by
applications in probabilistic forecasting, and thus we assess
combination formulas and aggregation methods in terms of coherency,
calibration, and dispersion.

While it is not difficult to construct coherent combination formulas,
the formulas typically used in practice tend to be incoherent.  Our
analytical results about the linear pool extend and unify extant work
by Dawid et al.~(1995), Hora (2004) and Ranjan and Gneiting (2010),
and show that any linear combination formula with strictly positive
coefficients fails to be coherent.  In this sense, combined evidence
is nonlinear.  An interesting open problem is whether or not
parsimonious ramifications, such as spread-adjusted or beta
transformed linear combination formulas, are coherent.

That said, the applied relevance of the notion of coherency is
limited, as predictive distributions or expert opinions issued in
practice are hardly ever ideal.  In typical practice, underdispersed
or approximately neutrally dispersed predictive distributions are to
be aggregated.  In the case of underdispersed components, the tendency
of linear combination formulas to increase dispersion can be
beneficial, and helps explain the success of linear pooling in
applications (Madigan and Raftery 1994).  However, if the components
are neutrally dispersed, the failure of the traditional linear pool to
be flexibly dispersive is a serious limitation.  Berrocal et
al.~(2007), Glahn et al.~(2009) and Kleiber et al.~(2011) observed
this empirically in the context of probabilistic weather forecasts,
and proposed a special case of the spread-adjusted linear pool as an
ad hoc remedy.  Our theoretical results document the increased
flexibility of the spread-adjusted linear pool, and demonstrate that
the beta-transformed linear pool is exchangeably flexibly dispersive.

Not surprisingly, the parsimonity principle and the bias-variance
tradeoff apply in the practice of the combination of predictive
distributions.  Thus, in data poor settings, where training data are
scarce, the parsimonious traditional linear pool might be the method
of choice, despite its theoretical shortcomings, as demonstrated
persuasively in the recent simulation study of Clements and Harvey
(2011).  In data rich settings, where predictive models can reliably
be estimated, linear aggregation tends to be suboptimal.  Hence, we
have studied parsimonious nonlinear alternatives, including the
spread-adjusted linear pool (SLP) and the beta transformed linear pool
(BLP).  Further options include consensus methods (Winkler 1968) and
nonparametric approaches, such as isotonic recursive partitioning
(Luss, Rosset and Shahar 2011).  As Winkler (1986) noted, ``different
combining rules are suitable for different situations''.

The SLP and BLP approaches can also be used to provide calibration and
dispersion corrections to a single predictive distribution, similar to
the methods described by Cox (1958), Platt (1999), Zadrozny and Elkan
(2002) and Primo et al.~(2009) in the context of probability forecasts
for a binary event.  An interesting question then is whether
dispersion adjustments ought to be applied to the individual
components prior to the aggregation.  In situations in which the
components show substantially differing degrees of dispersion, or are
uniformly under- or overdispersed, we indeed see potential benefits in
doing this, with (here unreported) simulation experiments providing
partial support to this view.  In our temperature example, the
components derive from point forecasts, which is the most extreme form
of underdispersion, and prior to aggregating the components we apply a
simple Gaussian technique, which obtains approximately neutrally
dispersed individual density forecasts.

In addition to their relevance in probabilistic forecasting, our
findings bear on the related problem of the fusion of expert opinions
that are expressed in the form of probability distributions.  Ha-Duong
(2008) reviews methods for doing this, and applies them to combine
expert opinions about the climate sensitivity constant, which is a key
quantity in the study of the greenhouse effect.  Our analytic results
imply that if each individual expert is neutrally dispersed, linear
aggregation results in combined assessments that are underconfident
and show an unduly wide range of uncertainty, when in fact a sharper
assessment could be made.  In situations of this type, the
parsimonious SLP technique with its single, easily interpretable
spread-adjustment parameter, $c$, might be preferable to the BLP,
particularly when the individual components are symmetric and
unimodal.  Our theoretical results show that neutrally dispersed
components require values of $c < 1$, with more specific
recommendations depending on the subject matter and situation at hand.
The case $c = 1$ corresponds to the traditional linear pool and can be
a useful choice in situations in which overconfident experts make
underdispersed judgements.

\section*{Appendix A: Details for Example \ref{ex:dispersion}}

Let $Z_\sigma = F_\sigma(Y)$ denote the probability integral transform
of the CDF-valued random quantity $F_\sigma$.  Then the random
variable $Z_\sigma$ has expectation $\frac{1}{2}$ and its cumulative
distribution function is $H_\sigma(z) = \Phi \hsp (\sigma \,
\Phi^{-1}(z))$.  In particular, $Z_1$ is uniformly distributed.  If
$\sigma < 1$ then $|Z_\sigma - \frac{1}{2}|$ is stochastically larger
than $|Z_1 - \frac{1}{2}|$ and therefore
\[
{
\textstyle
\var(Z_\sigma)
= \myE \hsp (Z_\sigma - \myE[Z_\sigma])^2
= \myE \hsp |Z_\sigma - \frac{1}{2}|^2
> \myE \hsp |Z_1 - \frac{1}{2}|^2
}
= \frac{1}{12}.
\]
An analogous argument applies when $\sigma > 1$.  To prove the
variance formula (\ref{eq:var}), we use the fact that $\var(Z_\sigma)
= \myE [Z_\sigma^2] - (\myE[Z_\sigma])^2$ and invoke the well-known
expectation equality $ \myE [Z^r] = r \int_0^\infty z^{r-1} (1 - H(z))
\: \dd z $ for a nonnegative random variable $Z$ with cumulative
distribution function $H$, where $r > 0$.

\section*{Appendix B: Method of scoring}

Here we give details for the method of scoring (see, for example,
Ferguson 1996) for numerically maximizing the mean logarithmic score
or log likelihood function (\ref{eq:log.likelihood}) of the BLP model
as a function of the nonnegative weights $w_1, \ldots, w_k$ that sum
to 1, and transformation parameters $\alpha, \beta > 0$.  Let $Y$
denote a random variable that has a beta distribution with parameters
$\alpha$ and $\beta$.  Then
\begin{eqnarray*}
\frac{\partial \ell}{\partial \alpha} & = &
   \sum_{j=1}^J \log \! \left( \sum_{i=1}^k w_i \hsp F_{ij}(y_j) \right)
   - \hsp J \, \myE \hsp [\log Y], \\
\frac{\partial \ell}{\partial \beta} & = &
   \sum_{j=1}^J \log \! \left( 1 - \sum_{i=1}^k w_i \hsp F_{ij}(y_j) \right)
    - \hsp J \, \myE \hsp [\log(1-Y)]
\end{eqnarray*}
and
\[
 \frac{\partial \ell}{\partial w_i}
= \sum_{j=1}^J \left(
  \frac{(\alpha-1)(F_{ij}(y_j) - F_{kj}(y_j))}{\sum_{l=1}^k w_l \hsp F_{lj}(y_j)}
- \frac{(\beta-1)(F_{ij}(y_j) - F_{kj}(y_j))}{1 - \sum_{l=1}^k w_l \hsp F_{lj}(y_j)}
+ \frac{f_{ij}(y_j) - f_{kj}(y_j)}{\sum_{l=1}^k w_l \hsp f_{lj}(y_j)} \right)
\]
for $i = 1, \ldots, k-1$.  The second derivatives are
\[
\frac{\partial^2 \ell}{\partial \alpha^2} = - \, J \, \var\hsp(\log(Y)),
\;\;
\frac{\partial^2 \ell}{\partial \beta^2} = - \, J \, \var \hsp (\log(1- Y)),
\;\;
\frac{\partial^2 \ell}{\partial \alpha \, \partial \beta} = - \, J \, \cov \hsp (\log(Y),\log(1-Y))
\]
and
\[
\frac{\partial^2 \ell}{\partial \alpha \, \partial w_i}
= \sum_{j=1}^J \frac{F_{ij}(y_j) - F_{kj}(y_j)}{\sum_{l=1}^k w_l \hsp F_{lj}(y_j)},
\qquad
\frac{\partial^2 \ell}{\partial \beta \, \partial w_i}
= \sum_{j=1}^J \frac{F_{kj}(y_j) - F_{ij}(y_j)}{1 - \sum_{l=1}^k w_l \hsp F_{lj}(y_j)}
\]
for $i = 1, \ldots, k-1$, while
\begin{eqnarray*}
\lefteqn{\hspace{-78.5mm}
\frac{\partial^2 \ell}{\partial w_{i_1} \partial w_{i_2}} \; = \;
- \sum_{j=1}^J \frac{(f_{{i_1}j}(y_j) - f_{kj}(y_j)) \hsp (f_{{i_2}j}(y_j) - f_{kj}(y_j))}
                    {(\sum_{l=1}^k w_l \hsp f_{lj}(y_j))^2}} \\
\lefteqn{\hspace{-78.5mm} - \sum_{j=1}^J \left(
  \frac{\alpha - 1}{(\sum_{l=1}^k w_l \hsp F_{lj}(y_j))^2} +
  \frac{\beta - 1}{(1- \sum_{l=1}^k w_l \hsp F_{lj}(y_j))^2} \right)
  \left( F_{{i_1}j}(y_j) - F_{kj}(y_j) \right)
  \left( F_{{i_2}j}(y_j) - F_{kj}(y_j) \right)}
\end{eqnarray*}
for $i_1 = 1, \ldots, k-1$ and $i_2 = 1, \ldots, k-1$.  The method of
scoring now applies Newton's algorithm to optimize the likelihood as a
function of the parameter vector.

\section*{Acknowledgements}

This research was supported by the Alfried Krupp von und zu Behlen
Foundation and by the United States National Science Foundation under
Awards ATM-0724721 and DMS-0706745 to the University of Washington.
Part of it was done at the University of Washington and at the Newton
Institute in Cambridge.  We are grateful to Cliff Mass and Jeff Baars
for providing the temperature data, and to Hajo Holzmann, Don
Percival, Adrian Raftery, Michael Scheuerer, Thordis Thorarinsdottir,
Alexander Tsyplakov, Bob Winkler and Johanna Ziegel for discussions
and preprints.  Of course, the usual disclaimer applies.

\section*{References}

\newenvironment{reflist}{\begin{list}{}{\itemsep 0mm \parsep 1mm
\listparindent -7mm \leftmargin 7mm} \item \ }{\end{list}}

\vspace{-7.5mm}
\begin{reflist}

Berrocal, V.~J., Raftery, A.~E.~and Gneiting, T.~(2007).  Combining
spatial statistical and ensemble information for probabilistic weather
forecasting.  {\em Monthly Weather Review}, {\bf 135}, 1386--1402.

Bj{\o}rnland, H.~C., Gerdrup, K., Jore, A.~S., Smith, C.~and Thorsrud,
L.~A.~(2011).  Weights and pools for a Norwegian density combination.
{\em North American Journal of Economics and Finance}, {\bf 22},
61--76.

Bollerslev, T.~(1986).  Generalized autoregressive conditional
heteroscedasticity, {\em Journal of Econometrics}, {\bf 31}, 307--327.

Br\"ocker, J., Siegert, S.~and Kantz, H.~(2011).  Comment on
``Conditional exceedance probabilities''.  {\em Monthly Weather
  Review}, {\bf 139}, in press.

Brockwell, A.~E.~(2007).  Universal residuals: A multivariate
transformation, {\em Statistics and Probability Letters}, {\bf 77},
1473--1478.

Clemen, R.~T.~and Winkler, R.~L.~(1999).  Combining probability
distributions from experts in risk analysis.  {\em Risk Analysis},
{\bf 19}, 187--203.

Clemen, R.~T.~and Winkler, R.~L.~(2007).  Aggregating probability
distributions.  In Ward, E., Miles, R.~F.~and von Winterfeldt,
D.~(eds.), {\em Advances in Decision Analysis: From Foundations to
  Applications}, Cambridge University Press, pp.~154--176.

Clements, M.~P.~and Harvey, D.~I.~(2011).  Combining probability
forecasts.  {\em International Journal of Forecasting}, {\bf 27},
208--223.

Corradi, V.~and Swanson, N.~R.~(2006).  Predictive density evaluation.
In Elliott, G., Granger, C.~W.~J.~and Timmermann, A.~(eds.), {\em
  Handbook of Economic Forecasting}, vol.~1, Elsevier North-Holland,
pp.~197--284.

Cox, D.~R.~(1958).  Two further applications of a model for binary
regression.  {\em Biometrika}, {\bf 45}, 562--565.

Czado, C., Gneiting, T.~and Held, L.~(2009).  Predictive model
assessment for count data.  {\em Biometrics}, {\bf 65}, 1254--1261.

Dawid, A.~P.~(1984).  Statistical theory: The prequential approach
(with discussion and rejoinder).  {\em Journal of the Royal
  Statistical Society Ser.~A}, {\bf 147}, 278--290.

Dawid, A.~P.~(1986).  Probability forecasting, in Kotz, S., Johnson,
N.~L.~and Read, C.~B.~(eds.), {\em Encyclopedia of Statistical
  Sciences}, Vol.~7, Wiley, New York, pp.~210--218.

Dawid, A.~P., DeGroot, M.~H.~and Mortera, J.~(1995).  Coherent
combination of experts' opinions (with discussion and rejoinder).
{\em Test}, {\bf 4}, 263--313.

Diebold, F.~X., Gunther, T.~A.~and Tay, A.~S.~(1998).  Evaluating
density forecasts with applications to financial risk management.
{\em International Economic Review}, {\bf 39}, 863--883.

Dutton, J.~A.~(2002).  Opportunities and priorities in a new era for
weather and climate services.  {\em Bulletin of the American
Meteorological Society}, {\bf 83}, 1303--1311.

Eckel, F.~A.~and Mass, C.~F.~(2005).  Aspects of effective short-range
ensemble forecasting.  {\em Weather and Forecasting}, {\bf 20},
328--350.

Ferguson, T.~S.~(1996).  {\em A Course in Large Sample Theory}.
Chapman and Hall.

Fraley, C., Raftery, A.~E., Gneiting, T., Berrocal, V.~J.~and
Sloughter, J.~M.~(2011).  Probabilistic weather forecasting in {\sc
  R}.  {\em R Journal}, in press.

French, S.~and R\'{i}os Insua, D.~(2000).  {\em Statistical Decision
  Theory}.  Arnold, London.

Garratt, A., Mitchell, J., Vahey, S.~P.~and Wakerly, E.~C.~(2011).
Real-time inflation forecast densities from ensemble Phillips curves.
{\em North American Journal of Economics and Finance}, {\bf 22},
77--87.

Genest, C.~(1984a).  A characterization theorem for externally
Bayesian groups.  {\em Annals of Statistics}, {\bf 12}, 1100--1105.

Genest, C.~(1984b).  A conflict between two axioms for combining
subjective distributions. {\em Journal of the Royal Statistical
  Society Ser.~B}, {\bf 46}, 403--405.

Genest, C.~and MacKay, R.~J.~(1986).  Copules archim\'ediennes et
families de lois bidimensionnelles dont les marges sont donn\'ees.
{\em Canadian Journal of Statistics}, {\bf 14}, 145--159.

Genest, C.~and Zidek, J.~(1986).  Combining probability distributions:
A critique and an annotated bibliography.  {\em Statistical Science},
{\bf 1}, 114--135.

Genest, C., McConway, K.~J.~and Schervish, M.~J.~(1986).
Characterization of externally Bayesian pooling operators.  {\em
  Annals of Statistics}, {\bf 14}, 487--501.

Gewewke, J.~and Amisano, G.~(2011).  Optimal prediction pools.  {\em
  Journal of Econometrics}, in press,
\url{doi:10.1016/j.jeconom.2011.02.017}.

Glahn, B., Peroutka, M., Wiedenfeld, J., Wagner, J., Zylstra, G.,
Schuknecht, B.~and Jackson, B.~(2009).  MOS uncertainty estimates in
an ensemble framework.  {\em Monthly Weather Review}, {\bf 137},
246--268.

Gneiting, T.~(2008).  Editorial: Probabilistic forecasting.  {\em
  Journal of the Royal Statistical Society Ser.~A}, {\bf 171},
319--321.

Gneiting, T.~and Raftery, A.~E.~(2005).  Weather forecasting with
ensemble methods.  {\em Science}, {\bf 310}, 248--249.

\newpage
Gneiting, T.~and Raftery, A.~E.~(2007).  Strictly proper scoring
rules, prediction, and estimation.  {\em Journal of the American
Statistical Association}, {\bf 102}, 359--378.

Gneiting, T., Balabdaoui, F.~and Raftery, A.~E.~(2007).  Probabilistic
forecasts, calibration and sharpness.  {\em Journal of the Royal
Statistical Society Ser.~B}, {\bf 69}, 243--268.

Ha-Duong, M.~(2008).  Hierarchical fusion of expert opinions in the
transferable belief model, application to climate sensitivity.  {\em
International Journal of Approximate Reasoning}, {\bf 49}, 555--574.

Hall, S.~G.~and Mitchell, J.~(2007).  Combining density forecasts.
{\em International Journal of Forecasting}, {\bf 23}, 1--13.

Hamill, T.~M.~(2001).  Interpretation of rank histograms for verifying
ensemble forecasts.  {\em Monthly Weather Review}, {\bf 129},
550--560.

Held, L., Rufibach, K.~and Balabdaoui, F.~(2010).  A score regression
approach to assess calibration of continuous probabilistic
predictions.  {\em Biometrics}, {\bf 66}, 1295--1305.

Hoeting, J.~A., Madigan, D., Raftery, A.~E.~and Volinsky, C.~T.~(1999)
Bayesian model averaging: A tutorial (with discussion).  {\em
  Statistal Science}, {\bf 4}, 382--417.

Hora, S.~C.~(2004).  Probability judgements for continuous quantities:
Linear combinations and calibration.  {\em Management Science}, {\bf
50}, 597--604.

Hora, S.~C.~(2010).  An analytic method for evaluating the performance
of aggregation rules for probability densities.  {\em Operations
  Research}, {\bf 58}, 1440--1449.

Jore, A.~S., Mitchell, J.~and Vahey, S.~P.~(2010).  Combining forecast
densities from VARs with uncertain stabilities.  {\em Journal of
  Applied Econometrics}, {\bf 25}, 621--634.

Jouini, M.~N.~and Clemen, R.~T.~(1996).  Copula models for aggregating
expert opinions.  {\em Operations Research}, {\bf 44}, 444--457.

Kascha, C.~and Ravazzolo, F.~(2010).  Combining inflation density
forecasts.  {\em Journal of Forecasting}, {\bf 29}, 231--250.

Kleiber, W., Raftery, A.~E., Baars, J., Gneiting, T., Mass, C.~F.~and
Grimit, E.~(2011).  Locally calibrated probabilistic temperature
forecasting using geostatistical model averaging and local Bayesian
model averaging.  {\em Monthly Weather Review}, in press.

Luss, R., Rosset, S.~and Shahar, M.~(2011).  Isotonic recursive
partitioning.  Preprint, \url{arXiv:1102.5496}.

Madigan, D.~and Raftery, A.~E.~(1994).  Model selection and accounting
for model uncertainty in graphical models using Occam's window.  {\em
  Journal of the American Statistical Association}, {\bf 89},
1535--1546.

Mason, S.~J., Galpin, S., Goddard, L., Graham, N.~E.~and Rajaratnam,
B.~(2007).  Conditional exceedance probabilities.  {\em Monthly
  Weather Review}, {\bf 135}, 363--372.

Mason, S.~J., Tippett, M.~K., Weigel, A.~P., Goddard, L.~and
Rajaratnam, B.~(2011).  Reply to ``Comment on `Conditional exceedance
probabilities'{''}.  {\em Monthly Weather Review}, {\bf 139}, in press.

\newpage
Matheson, J.~E.~and Winkler, R.~L.~(1976).  Scoring rules for
continuous probability distributions.  {\em Management Science}, {\bf
22}, 1087--1096.

McConway, K.~J.~(1981).  Marginalization and linear opinion pools.
{\em Journal of the American Statistical Association}, {\bf 76},
410--414.

McNeil, A.~J.~and Ne\u{s}lehov\'a, J.~(2009).  Multivariate
Archimedean copulas, $d$-monotone functions and $\ell_1$-symmetric
functions.  {\em Annals of Statistics}, {\bf 37}, 3059--3097. 

Mitchell, J.~and Hall, S.~G.~(2005).  Evaluating, comparing and
combining density forecasts using the KLIC with an application to the
Bank of England and NIESR `fan' charts of inflation.  {\em Oxford
Bulletin of Economics and Statistics}, {\bf 67}, 995--1033.

Mitchell, J.~and Wallis, K.~F.~(2010).  Evaluating density forecasts:
forecast combinations, model mixtures, calibration and sharpness. 
{\em Journal of Applied Econometrics}, \url{DOI:10.1002/jae.1192}.

Murphy, A.~H.~and Winkler, R.~L.~(1987).  A general framework for
forecast verification.  {\em Monthly Weather Review}, {\bf 115},
1330--1338.

Murphy, A.~H.~and Winkler, R.~L.~(1992).  Diagnostic verification of
probability forecasts.  {\em International Journal of Forecasting},
{\bf 7}, 435--455.

Palmer, T.~N.~(2002).  The economic value of ensemble forecasts as a
tool for risk assessment: From days to decades.  {\em Quarterly
Journal of the Royal Meteorological Society}, {\bf 128}, 747--774.

Platt, J.~C.~(1999).  Probabilistic outputs for support vector
machines and comparisons to regularized likelihood methods.  In Smola,
A., Bartlett, P., Sch\"olkopf, B.~and Schuurmans, D., eds., {\em
  Advances in Large Margins Classifiers}, MIT Press.

Primo, C., Ferro, C.~A.~T., Jolliffe, I.~T.~and Stephenson,
D.~B.~(2009).  Calibration of probabilistic forecasts of binary
events.  {\em Monthly Weather Review}, {\bf 137}, 1142--1149.

Raftery, A.~E., Gneiting, T., Balabdaoui, F.~and
Polakowski,~M.~(2005).  Using Bayesian model averaging to calibrate
forecast ensembles.  {\em Monthly Weather Review}, {\bf 133},
1155--1174.

Ranjan, R.~and Gneiting, T.~(2010).  Combining probability forecasts.
{\em Journal of the Royal Statistical Society Ser.~B}, {\bf 72},
71--91.

R Development Core Team (2011).  R: A language and environment for statistical
computing.  R Foundation for Statistical Computing, Vienna, Austria,
ISBN 3-900051-07-0, \url{http://www.R-project.org}.

Rosenblatt, M.~(1952).  Remarks on a multivariate transformation.
{\em Annals of Mathematical Statistics}, {\bf 23}, 470--472.

Schervish, M.~J.~(1989).  A general method for comparing probability
assessors.  {\sl Annals of Statistics}, {\bf 17}, 1856--1879.

Sloughter, J.~M., Raftery, A.~E., Gneiting, T.~and Fraley, C.~(2007).
Probabilistic quantitative precipitation forecasting using Bayesian
model averaging.  {\em Monthly Weather Review}, {\bf 135}, 3209--3220.

Stone, M.~(1961).  The linear pool.  {\em Annals of Mathematical
Statistics}, {\bf 32}, 1339--1342.

Tsyplakov, A.~(2011).  Evaluating density forecasts: a comment.  MPRA
paper no.~31233, \url{http://mpra.ub.uni-muenchen.de/31233}.

Wallis, K.~F.~(2005).  Combining density and interval forecasts: A
modest proposal.  {\em Oxford Bulletin of Economics and Statistics},
{\bf 67}, 983--994.

Winkler, R.~L.~(1968).  The consensus of subjective probability
distributions.  {\em Management Science}, {\bf 15}, B61--B75.

Winkler, R.~L.~(1981).  Combining probability distributions from
dependent information sources.  {\em Management Science}, {\bf 27},
479--488.

Winkler, R.~L.~(1986).  Comment on ``Combining probability
distributions: A critique and an annotated bibliography''.  {\em
  Statistical Science}, {\bf 1}, 138--140.

Wuertz, D.~and Rmetrics Core Team (2007). The fGarch Package.
Reference manual, available at
\url{http://www.mirrorservice.org/sites/lib.stat.cmu.edu/R/CRAN/doc/packages/fGarch.pdf}.

Zadrozny, B.~and Elkan, C.~(2002).  Transforming classifier scores
into accurate multiclass probability estimates.  Proceedings of the
eighth ACM SIGKDD International Conference on Knowledge Discovery and
Data Mining, pp.~694--699. 

Zarnowitz, V.~(1969).  The new ASA-NBER survey of forecasts by
economic statisticians.  {\em American Statistician}, {\bf 23},
12--16.
\end{reflist}

\end{document}